\input amstex
\documentstyle{amsppt}
\loadbold

\magnification=\magstep1

\TagsAsMath
\TagsOnLeft
\NoBlackBoxes


\pageheight{9.0truein}
\pagewidth{6.5truein}

\long\def\ignore#1\endignore{#1}

\ignore \input xy \xyoption{matrix} \xyoption{arrow} \endignore

\def\edge{\ar@{-}}
\def\place{*{}+0}
\def\plb{*{\bullet}+0}
\def\plc{*{\circ}+0}

\def\hrz{\place \edge[r] &\place}

\def\hrzvrt{\place \edge[r] \edge[d] &\place \edge[d]}

\def\dropup#1#2{\save+<0ex,#1ex> \drop{#2} \restore}
\def\dropleft#1#2{\save+<-#1ex,0ex> \drop{#2} \restore}

\def\B{{\Cal B}}
\def\C{{\Cal C}}
\def\E{{\Cal E}}

\def\O{{\Cal O}}

\def\CC{{\Bbb C}}
\def\NN{{\Bbb N}}
\def\ZZ{{\Bbb Z}}
\def\qhat{\widehat{q}}
\def\bfc{{\boldsymbol c}}
\def\bfr{{\boldsymbol r}}

\def\bftau{{\boldsymbol \tau}}

\def\dtil{\widetilde{d}}

\def\util{\widetilde{u}}

\def\wtil{\widetilde{w}}
\def\Htil{\widetilde{H}}
\def\Itil{\widetilde{I}}
\def\Jtil{\widetilde{J}}
\def\Ptil{\widetilde{P}}

\def\Xtil{\widetilde{X}}

\def\Xbar{\overline{X}}

\def\Oq{{\Cal O}_q}
\def\OqMn{\Oq(M_n(k))}

\def\OqMtwo{\Oq(M_2(k))}
\def\OqMthree{\Oq(M_3(k))}

\def\OqMUV{\Oq(M_{U,V}(k))}
\def\OqSLn{\Oq(SL_n(k))}
\def\OqGLn{\Oq(GL_n(k))}
\def\OqG{\Oq(G)}
\def\GLn{GL_n(k)}
\def\pmqdot{\pm q^{\bullet}}
\def\kx{k^\times}
\def\rc{\bfr,\bfc}
\def\RC{{\boldsymbol R}{\boldsymbol C}}
\def\Krc{K_{\rc}}
\def\Arc{A_{\rc}}
\def\Drc{D_{\rc}}
\def\Drctil{{\widetilde D}_{\rc}}
\def\drc#1{d^{\rc}_{#1}}
\def\drctil#1{\dtil^{\rc}_{#1}}
\def\drctilinv#1{(\dtil^{\rc}_{#1})^{-1}}

\def\onen{\{1,\dots,n\}}
\def\Rrpluszero{R^+_{\bfr,0}}
\def\Rcminuszero{R^-_{\bfc,0}}
\def\Rrplus{R^+_{\bfr}}
\def\Rcminus{R^-_{\bfc}}
\def\ror{\Rrplus \otimes \Rcminus}
\def\pirpluszero{\pi^+_{\bfr,0}}
\def\picminuszero{\pi^-_{\bfc,0}}
\def\betarc{\beta_{\rc}}
\def\Brc{B_{\rc}}
\def\betatilrc{\widetilde{\beta}_{\rc}}
\def\Brplus{B^+_{\bfr}}
\def\Bcminus{B^-_{\bfc}}

\def\OqMtwothreemodboxbox{\Oq(M_{2,3}(k))/ \langle \square
   \kern-.1667em \square \rangle}

\def\congmodKrc{\mathrel{\equiv_{\rc}}}
\def\Rt{{\boldsymbol R}_t}

\def\fract{\operatorname{Fract}}
\def\spec{\operatorname{spec}}
\def\specrc{\spec_{\rc}}

\def\Hspec{H\text{-}\operatorname{spec}}

\def\Aut{\operatorname{Aut}}
\def\medcap{\operatornamewithlimits{\tsize{\bigcap}}}

\def\barcelona{{\bf 1}}
\def\Cau{{\bf 2}}
\def\DeCL{{\bf 3}}
\def\DeCP{{\bf 4}}
\def\GLmurcia{{\bf 5}}
\def\Duke{{\bf 6}}
\def\GLenseq{{\bf 7}}
\def\GLet{{\bf 8}}
\def\specstrat{{\bf 9}}
\def\HoLeone{{\bf 10}}
\def\HoLetwo{{\bf 11}}
\def\Josbook{{\bf 12}}
\def\LeSt{{\bf 13}}
\def\MoRe{{\bf 14}}
\def\NYM{{\bf 15}}
\def\PaWa{{\bf 16}}
\def\Von{{\bf 17}}

\topmatter

\title Prime ideals invariant under winding automorphisms in quantum
matrices\endtitle

\rightheadtext{Winding-in\-var\-i\-ant prime ideals in quantum matrices}

\author K. R. Goodearl and T. H. Lenagan\endauthor

\address Department of Mathematics, University of California, Santa
Barbara, CA 93106, USA\endaddress
\email goodearl\@math.ucsb.edu\endemail

\address Department of Mathematics, J.C.M.B., Kings Buildings, Mayfield
Road, Edinburgh EH9 3JZ, Scotland\endaddress
\email tom\@maths.ed.ac.uk\endemail

\thanks This research was partially supported by NSF
research grants DMS-9622876 and DMS-9970159 and by NATO Collaborative
Research Grant CRG.960250. Some of the research was done while both
authors were visiting the Mathematical Sciences Research Institute in
Berkeley during the winter of 2000, and they thank MSRI for its
support.
\endthanks

\abstract The main goal of the paper is to establish the existence of
tensor product decompositions for those prime ideals $P$ of the algebra
$A=\OqMn$ of quantum $n\times n$ matrices which are invariant under
winding automorphisms of $A$, in the generic case ($q$ not a root of
unity). More specifically, every such $P$ is the kernel of a map of the
form
$$A \longrightarrow A\otimes A \longrightarrow A^+\otimes A^-
\longrightarrow (A^+/P^+)\otimes (A^-/P^-)$$
where $A\rightarrow A\otimes A$ is the comultiplication, $A^+$ and $A^-$
are suitable localized factor algebras of $A$, and $P^\pm$ is a prime
ideal of $A^\pm$ invariant under winding automorphisms. Further, the
algebras $A^\pm$, which vary with $P$, can be chosen so that the
correspondence $(P^+,P^-) \mapsto P$ is a bijection. The main theorem is
applied, in a sequel to this paper, to completely determine the
winding-in\-var\-i\-ant prime ideals in the generic quantum $3\times3$ matrix
algebra. 
\endabstract

\endtopmatter

\document

\head Introduction \endhead

This paper represents part of an ongoing project to determine the prime
and primitive spectra of the generic quantized coordinate ring of
$n\times n$ matrices, $\OqMn$. Here $k$ is an arbitrary field and
$q\in\kx$ is a non-root of unity. The current intermediate goal is to
determine the prime ideals of $\OqMn$ invariant under all winding
automorphisms. (See below for a discussion of the relations between these
{\it winding-in\-var\-i\-ant\/} primes and the full prime spectrum of
$\OqMn$.) Our main result exhibits a bijection between these primes and
pairs of winding-in\-var\-i\-ant primes from certain `localized
step-triangular factors' of $\OqMn$, namely the algebras
$$\aligned \Rrplus &= \bigl( \OqMn/\langle X_{ij}\mid j>t \text{\ or\ }
i<r_j \rangle \bigr) [ \Xbar^{-1}_{r_11}, \dots, \Xbar^{-1}_{r_tt}] \\
\Rcminus &= \bigl( \OqMn/\langle X_{ij}\mid i>t \text{\ or\ }
j<c_i \rangle \bigr) [\Xbar^{-1}_{1c_1}, \dots, \Xbar^{-1}_{tc_t}]
\endaligned$$
where $\bfr= (r_1,\dots,r_t)$ and $\bfc=
(c_1,\dots,c_t)$ are strictly increasing sequences of integers in the
range $1,2,\dots,n$. In particular, since each $\Rrplus$ and $\Rcminus$
can be presented as a skew-Laurent extension of a localized factor
algebra of
$\Oq(M_{n-1}(k))$, the above bijection can be used to obtain
descriptions (as pullbacks of primes in the algebras $\ror$) of the
winding-in\-var\-i\-ant primes of
$\OqMn$ from those of
$\Oq(M_{n-1}(k))$. In a sequel \cite{\GLenseq} to this paper, we follow
the route just sketched to develop a complete list, with sets of
generators, of the winding-in\-var\-i\-ant primes in $\OqMthree$.

The theorem indicated above depends on some detailed structural results
concerning $\OqMn$ and on some general work with primes in tensor
product algebras invariant under group actions. First, we construct a
partition of
$\spec \OqMn$ indexed by pairs $(\rc)$ as above, together with
localized factor algebras $\Arc$ of $\OqMn$, such that the portion of
$\spec\OqMn$ indexed by $(\rc)$ is Zariski-homeomorphic to $\spec\Arc$.
We next prove that $\Arc$ is isomorphic to a subalgebra
$\Brc$ of
$\ror$, identify the structure of $\Brc$, and show that $\ror$ is a
skew-Laurent extension of $\Brc$. Finally, with the help of some
general work on tensor products, we prove that each winding-in\-var\-i\-ant
prime of $\Brc$ extends uniquely to a winding-in\-var\-i\-ant prime of
$\ror$, and that the latter primes can be uniquely expressed in the
form $(P^+\otimes \Rcminus) +(\Rrplus\otimes P^-)$ where $P^+$
(respectively, $P^-$) is a winding-in\-var\-i\-ant prime in $\Rrplus$
(respectively, $\Rcminus$). We thus conclude that every winding-in\-var\-i\-ant
prime of $\OqMn$ can be uniquely expressed as the kernel of a map
$$\OqMn \longrightarrow \OqMn\otimes\OqMn \longrightarrow \ror
\longrightarrow (\Rrplus/P^+) \otimes (\Rcminus/P^-),$$
where the first arrow is comultiplication and the others are tensor
products of localization or quotient maps.

\definition{Algebraic background} The algebra $\OqMn$ has standard
generators
$X_{ij}$ for $i,j=1,\dots,n$ and relations which we recall in
(5.1)(a), along with the bialgebra structure of this algebra. The latter
structure allows us to define {\it left and right winding
automorphisms\/} corresponding to those characters ($k$-algebra
homomorphisms $\OqMn
\rightarrow k$) which are invertible in $\OqMn^*$ with respect to the
convolution product (cf\. \cite{\barcelona, (I.9.25)} or \cite{\Josbook,
(1.3.4)} for the Hopf algebra case). It is well known that the collection
of left (respectively, right) winding automorphisms of $\OqMn$ forms a
group isomorphic to the diagonal subgroup of $GL_n(k)$, whose action on
the matrix of generators
$(X_{ij})$ is given by left (respectively, right) multiplication.
We combine these actions to obtain an action of the group $H= (\kx)^n
\times (\kx)^n$ on $\OqMn$ by $k$-algebra automorphisms satisfying the
rule
$$(\alpha_1,\dots,\alpha_n,\beta_1,\dots,\beta_n){.}X_{ij}=
\alpha_i\beta_jX_{ij}. \tag*$$
One indication of the extent of the symmetry given by this action is
the fact that there are only finitely many (actually, at most
$2^{n^2}$) primes of $\OqMn$ invariant under $H$ \cite{\specstrat,
(5.7)(i)}. The quoted result also shows that all $H$-primes of this
algebra are prime, and so the $H$-primes coincide with the
winding-invariant primes in $\OqMn$. (Recall that the definition of an
{\it $H$-prime\/} ideal is obtained from the standard ideal-theoretic
definition of a prime ideal by restricting to $H$-invariant ideals.)

In \cite{\specstrat, Theorem 6.6}, Letzter and the first author showed
that the overall picture of the prime spectrum of an algebra with
certain basic features like those of $\OqMn$ is determined to a great
extent by the primes invariant under a suitable group action. We quote
the improved version of this picture presented in \cite{\barcelona,
Theorem II.2.13}. Let $A$ be a noetherian algebra over an infinite
field $k$, and let $H=(\kx)^r$ be (the group of $k$-points of) an
algebraic torus acting rationally on $A$ by $k$-algebra automorphisms.
Each $H$-prime of $A$ is a prime ideal, and $\spec A$ is the
disjoint union of the sets
$$\spec_JA := \{ P\in \spec A \mid \medcap_{h\in H} h(P)= J\}$$
as $J$ ranges over the $H$-primes of $A$. Further:
\roster
\item"(1)" Let $\E_J$ denote the set of all regular $H$-eigen\-vec\-tors in
$A/J$. Then $\E_J$ is a denominator set, and the localization $A_J=
(A/J)[\E_J^{-1}]$ is an $H$-simple ring.
\item"(2)" $\spec_J A$ is homeomorphic to $\spec A_J$ via localization
and contraction.
\item"(3)" $\spec A_J$ is homeomorphic to $\spec Z(A_J)$ via
contraction and extension.
\item"(4)" $Z(A_J)$ is a Laurent polynomial ring, in at most $r$
indeterminates, over the fixed field $Z(\fract A/J)^H$.
\endroster
Under some additional hypotheses, satisfied by $\OqMn$, we also have:
\roster
\item"(5)" The primitive ideals of $A$ are exactly the maximal elements
of the sets $\spec_JA$.
\item"(6)" If $k$ is algebraically closed, then the primitive ideals
within each $\spec_JA$ are permuted transitively by $H$.
\endroster
Statement (5) was proved in \cite{\specstrat, Corollary 6.9} (see also
\cite{\barcelona, Theorem II.8.4}). Statement (6) is a consequence of
general transitivity theorems for algebraic group actions due to
Moeglin-Rentschler \cite{\MoRe, Th\'eor\`eme 2.12(ii)} and Vonessen
\cite{\Von, Theorem 2.2}, but the case where the acting group is a
torus is much easier (see \cite{\specstrat, Theorem 6.8} or
\cite{\barcelona, Theorem II.8.14}).

The above results indicate that to draw a complete picture of
$\spec\OqMn$, we need to determine the $H$-primes. That is easy to do in
case $n=2$; the result is recorded, for instance, in \cite{\GLmurcia,
(3.6)} (see \cite{\barcelona, Example II.2.14(d)} for more detail). In
general, we make the following
\enddefinition

\proclaim{Conjecture} Every $H$-prime of $\OqMn$ can be
generated by a set of quantum minors. \endproclaim

\noindent This conjecture is easily checked in case $n=2$ using the
information above, and we verify it for the case $n=3$ in
\cite{\GLenseq}. Further supporting evidence is
provided by recent work of Cauchon, who showed that distinct, comparable
$H$-primes in $\OqMn$ can be distinguished by the quantum minors they
contain \cite{\Cau, Proposition 6.2.2 and Th\'eor\`eme 6.2.1}. Another
source of support for the conjecture is the work of Hodges and Levasseur
\cite{\HoLeone,
\HoLetwo}, from which one can deduce that, up to certain localizations,
the winding-in\-var\-i\-ant primes of $\OqSLn$ are generated by quantum
minors (cf\.
\cite{\barcelona, Corollary II.4.12}). Since $\OqGLn$ is isomorphic to a
Laurent polynomial ring over $\OqSLn$
\cite{\LeSt, Proposition}, the above statement also holds in $\OqGLn$. In
particular, those $H$-primes of $\OqMn$ which do not contain the quantum
determinant can be generated, up to suitable localizations, by quantum
minors.

\definition{Geometric background} The classical origins of our main
theorem, especially as it applies to winding-invariant primes of $\OqMn$
not containing the quantum determinant, lie in the geometry of
`LU-decompositions' of invertible matrices. For this part of the
introduction, let us assume (to avoid complications) that $k$ is
algebraically closed. An LU-decomposition of a matrix $X\in \GLn$ is any
expression $X=LU$ where $L$ (respectively, $U$) is a lower
(respectively, upper) triangular invertible matrix. It is well known
that $X$ has such a decomposition if and only if the principal minors of
$X$ (those indexed by rows and columns from an initial segment of
$\{1,\dots,n\}$) are all nonzero. The LU-decomposable matrices 
thus form a dense open subvariety of $\GLn$, known as the {\it big
cell\/}. We may write the big cell in the form $B^+B^-$ where $B^+$
(respectively,
$B^-$) is the subgroup of lower (respectively, upper) triangular matrices
in $\GLn$, and we have
$$\O(B^+B^-)= \O(\GLn)[D^{-1}]$$
where $D$ is the multiplicative subset of $\O(\GLn)$ generated by the
principal minors. The comorphism of the multiplication map $B^+\times B^-
\rightarrow \GLn$ provides an embedding
$$\beta: \O(\GLn)[D^{-1}] \longrightarrow \O(B^+) \otimes \O(B^-);$$
the restriction of $\beta$ to $\O(\GLn)$ is just the composition of the
comultiplication map $\O(\GLn) \rightarrow \O(\GLn) \otimes \O(\GLn)$
with the tensor product of the restriction maps $\O(\GLn) \rightarrow
\O(B^\pm)$. The structure of the image of $\beta$ is easy to determine,
since $B^+\cap B^-$ is the diagonal subgroup of $\GLn$ and the subgroups
$B^\pm$ are semidirect products of their unipotent subgroups with this
diagonal subgroup. Namely, the image of $\beta$ is the subalgebra of
$\O(B^+) \otimes \O(B^-)$ generated by (the cosets of) the elements
$$\xalignat3 X_{ij}X_{jj}^{-1}\otimes1 &&1\otimes X_{ij}X_{ii}^{-1}
&&(X_{ii}\otimes X_{ii})^{\pm1}. & \tag\dagger \endxalignat$$
Further, $\O(B^+) \otimes \O(B^-)$ is a Laurent polynomial ring over the
image of $\beta$ with respect to indeterminates $1\otimes X_{ii}^{\pm1}$.

Quantum analogs of the above facts are known, but we have not been able
to locate complete statements in the literature. To formulate them, let
us write
$\Oq(B^+)$ and $\Oq(B^-)$ for the respective quotients of $\OqGLn$
modulo the ideals
$\langle X_{ij}
\mid i<j \rangle$ and $\langle X_{ij} \mid i>j \rangle$. Then:
\roster
\item"(a)" The
composition of the comultiplication map on $\OqGLn$ with
the tensor product of the quotient maps $\OqGLn \rightarrow \Oq(B^\pm)$
yields an embedding $\beta : \OqGLn
\rightarrow \Oq(B^+) \otimes \Oq(B^-)$ \cite{\PaWa, Theorem 8.1.1}.
\item"(b)" The
multiplicative subset $D$ of $\OqGLn$ generated by the principal quantum
minors is a denominator set.
\item"(c)" $\beta$ extends to $\OqGLn[D^{-1}]$, and the
image of this extension is the subalgebra $\B$ of $\Oq(B^+)
\otimes \Oq(B^-)$ generated by (the cosets of) the elements $(\dagger)$.
\item"(d)" $\Oq(B^+) \otimes \Oq(B^-)$ is a skew-Laurent extension of
$\B$ with respect to the variables $1\otimes X_{ii}^{\pm1}$.
\endroster
These facts will be proved as part of the case $\bfr=\bfc= (1,\dots,n)$
of our work below. To obtain them via existing results in the literature,
one first transfers the problem to
$\OqSLn$ using the isomorphism $\OqGLn \cong \OqSLn[z^{\pm1}]$
established in \cite{\LeSt, Proposition}; the desired conditions then
hold in the generality of $\OqG$, where $G$ is an arbitrary semisimple
algebraic group. The isomorphism of the appropriate localization of
$\OqG$ with the analog of $\B$ is given in \cite{\DeCL, Theorem 4.6} and 
\cite{\Josbook, Proposition 9.2.14}. It is easy to see that $\Oq(B^+)
\otimes \Oq(B^-)$ is a skew-Laurent extension of $\B$ with appropriate
variables; this is mentioned for the case where $q$ is a root of unity in
\cite{\DeCP, (4.6)}.

The above facts concerning $\OqGLn$ immediately carry over to $\OqMn$,
since $D$ is also a denominator set in that algebra and $\OqMn[D^{-1}]=
\OqGLn[D^{-1}]$. Whereas in $\OqGLn$ all prime ideals are disjoint from
$D$, that no longer holds in $\OqMn$, and a sequence of modified versions
of (a)--(d) are needed to yield information about those $H$-primes
of $\OqMn$ that meet $D$. These modifications, involving maps
$$\betarc: \OqMn\longrightarrow \OqMn\otimes\OqMn \longrightarrow \ror,$$
concern quantum analogs of what can be viewed as LU-decompositions for
certain locally closed subsets of $M_n(k)$ of the form $B^+w^+ew^-B^-$
where $w^\pm$ are permutation matrices and $e= e_{11}+\dots+e_{tt}$ is a
diagonal idempotent matrix. We leave the formulations of these geometric
facts to the interested reader.
\enddefinition

\definition{Some notation and conventions} Throughout the paper, let
$A=\OqMn$, where we fix a base field $k$, a positive integer $n$,
and a nonzero scalar
$q\in\kx$. See (5.1)(a) for the basic relations satisfied by the
standard generators $X_{ij}$ of $\OqMn$. As mentioned in (5.2), our
relations for $A$ differ from those in \cite{\PaWa} by an interchange of
$q$ and $q^{-1}$. On the other hand, they agree with the relations used
in \cite{\Cau, \NYM}. To match the relations in \cite{\HoLeone,
\HoLetwo}, replace $q$ by $q^2$.

While our main interest is in
the case that $q$ is not a root of unity, some of our results do not
require this assumption, and others require only that $q\ne\pm1$. Thus,
we impose no blanket hypotheses on
$q$. To simplify some formulas, write $\qhat= q-q^{-1}$. We also fix the
torus $H= (\kx)^n \times
(\kx)^n$ and its action on $A$ by winding automorphisms as described
in $(*)$ above. The algebra
$A$ is graded in a natural way by $\ZZ^{2n}= \ZZ^n\times\ZZ^n$, each
generator $X_{ij}$ having degree $(\epsilon_i,\epsilon_j)$ where
$\epsilon_1,\dots,\epsilon_n$ is the standard basis for $\ZZ^n$. We refer
to this grading as the {\it standard grading\/} on $A$. (As long as
$k$ is infinite, the homogeneous components of $A$ for this grading
coincide with the eigenspaces for the action of $H$.)

As in \cite{\Duke}, we use the notations $[I\mid J]$ and $[i_1\cdots
i_s\mid j_1\cdots j_s]$ for quantum minors in $\OqMn$, where $I$
(respectively, $i_1,\dots,i_s$) records the set (respectively, a list) of
the corresponding row indices, and similarly for column indices. Recall
that $[I\mid J]$ corresponds to the quantum determinant in a subalgebra
of $\OqMn$ isomorphic to $\Oq(M_s(k))$; we write $\Oq(M_{I,J}(k))$ for
that subalgebra (see (5.1)(c)). We allow the index sets $I$ and $J$ to be
empty, following the convention that $[\varnothing\mid \varnothing] =1$.

The symbols $\subset$ and $\subseteq$ will be reserved for proper and
arbitrary inclusions, respectively. We write $\sqcup$ to denote a disjoint
union.
\enddefinition

\head 1. A partition of $\spec A$ \endhead

We begin by investigating sets $\spec_{\rc}A$ of prime ideals defined by
certain `stepwise patterns' of quantum minors corresponding to strictly
increasing sequences $\bfr$ and $\bfc$ of row and column indices. Our
aim in this section is to show that these sets partition
$\spec A$ when $q$ is not a root of unity. In fact, as long as
$q\ne\pm1$, these sets at least partition the collection of completely
prime ideals of
$A$. (Recall from \cite{\GLet, Theorem 3.2} that when $q$ is not a root
of unity, all primes of $A$ are completely prime.) We proceed without
imposing special hypotheses on $q$ until needed.

\definition{1.1} We introduce the following partial ordering $\le$ on
index sets $I,I'\subseteq
\onen$ of the same cardinality. Write $I= \{i_1< \cdots<
i_l\}$ and $I'= \{i'_1< \cdots<
i'_l\}$; then $I\le I'$ if and only if $i_s\le i'_s$ for
$s=1,\dots\l$. (This is the same as the `column ordering' $\le_c$
used in \cite{\Duke}.) All order relations among index sets in this
paper will refer to the above partial ordering. This includes
statements that an index set with a particular property is minimal
among index sets with the same cardinality satisfying that property
(e.g., the sets $\Itil$ and $\Jtil$ in Theorem 1.9).
\enddefinition

\definition{1.2} Let $\RC$ denote the set of all pairs $(\rc)$ where
$\bfr$ and $\bfc$ are strictly increasing sequences in $\onen$ of the
same length, that is, $\bfr= (r_1,\dots,r_t)$ and $\bfc=
(c_1,\dots,c_t)$ in $\NN^t$ with $1\le r_1< r_2< \cdots< r_t\le n$
and $1\le c_1< c_2< \cdots< c_t\le n$. We allow $t=0$, in which case
$\bfr$ and $\bfc$ are empty sequences, denoted either $(\quad)$ or
$\varnothing$. When referring to the length $t$ of
$\bfr$ and $\bfc$, we write $(\rc) \in \RC_t$.

For $(\rc) \in \RC_t$, let $\Krc$ be the ideal of $A$ generated
by the following set of quantum minors:
$$\align \bigl\{ [I\mid J] \bigm| |I|>t \bigr\} &\cup \bigl\{
[I\mid J] \bigm| |I|=l\le t \text{\ and\ } I \not\ge
\{r_1,\dots,r_l\}
\bigr\} \\
 &\cup \bigl\{ [I\mid J] \bigm| |I|=l\le
t \text{\ and\ } J \not\ge \{c_1,\dots,c_l\} \bigr\}. \endalign$$
In particular, $K_{\varnothing,\varnothing}= \langle X_{ij} \mid
i,j\in \onen \rangle$.

Now set $\drc{l}= [r_1\cdots r_l \mid c_1\cdots c_l]$ for $l\le t$,
and observe that these quantum minors commute with each other (cf\.
(5.1)(c)). Let
$\Drc$ denote the multiplicative subset of $A$ generated by $\drc{1},
\dots, \drc{t}$. (In particular, $D_{\varnothing,\varnothing} =
\{1\}$.) Set $\drctil{l}= \drc{l}+\Krc \in A/\Krc$ for $l\le t$, and
let $\Drctil$ denote the image of $\Drc$ in $A/\Krc$.

We prove in this section that $\Drctil$ is a denominator set in
$A/\Krc$, and that when $q$ is generic, $\spec A$ is partitioned by
the subsets
$$\align \specrc A &:= \{ P\in \spec A \mid \Krc \subseteq P \text{\
and\ } P\cap \Drc= \varnothing \}\\
 &\hphantom{:}\approx \spec {(A/\Krc)[\Drctil^{-1}]} \endalign$$
as $(\rc)$ ranges over $\RC$.
\enddefinition

\proclaim{1.3\. Lemma} Let $I,J\subseteq \{1,\dots,n\}$ with
$|I|=|J|$, and set
$$L= \bigl\langle\; [I'\mid J'] \bigm| |I'|=|I|, \operatorname{\ and\ }
I'<I \operatorname{\ or\ } J'<J \;\bigr\rangle.$$
If $r,c\in \{1,\dots,n\}$ with $r\le \max(I)$ or $c\le \max(J)$, then
$$[I\mid J]X_{rc} -q^{2-\delta(r,I)-\delta(c,J)} X_{rc} [I\mid J] \ \in
\ L. \tag*$$
\endproclaim

\demo{Proof} If $r\in I$ and $c\in J$, then $[I\mid J]$ commutes with
$X_{rc}$ (cf\. Lemma 5.2(a)), and so (*) holds.

Next, suppose that $r\in I$ and $c\notin J$. If $j\in J$ and $j>c$,
then $J\sqcup \{c\} \setminus \{j\} <J$, whence $[I\mid J\sqcup \{c\}
\setminus \{j\}] \in L$ by definition of $L$. It follows from Lemma
5.2(b2) that $[I\mid J]X_{rc} -qX_{rc}
[I\mid J] \in L$. Thus (*) holds in this case. The case where $r\notin
I$ and $c\in J$ is proved similarly, using Lemma 5.2(c2).

Finally, suppose that $r\notin I$ and $c\notin J$, and note that either
$r< \max(I)$ or $c<\max(J)$. (This part of the proof is similar to that
of Lemma 5.7.) 
 There is no loss of generality in assuming that
$$I\sqcup \{r\}= J\sqcup \{c\}= \{1,\dots,n\},$$
whence either $r<n$ or $c<n$. In the notation of \cite{\PaWa, (4.3)},
$[I\mid J]= A(r\,c)$. Note that $\{1,\dots,n\} \setminus \{i\} <I$
when $i>r$, while $\{1,\dots,n\} \setminus \{j\} <J$ when $j>c$.
Hence, $A(i\,j)
\in L$ whenever either $i>r$ or $j>c$. 
Set $D_q= [1\cdots n\mid 1\cdots n]$. The basic q-Laplace relations
(Corollary 5.5) imply that $D_q$ lies in the ideal generated by all
the $A(n\,j)$, and in the ideal generated by all the $A(i\,n)$. Since
either $n>r$ or $n>c$, it follows that $D_q\in L$.

We now use the basic q-Laplace relations in the form given in
\cite{\PaWa, Corollary 4.4.4}. The first two relations yield
$$\sum_{j=1}^n (-q)^{j-r} X_{rj} A(r\,j)= \sum_{j=1}^n (-q)^{r-j}
A(r\,j) X_{rj} =D_q \in L. \tag1$$
Since $A(r\,j) \in L$ for $j>c$, we obtain the following congruences,
after multiplying the two sums in (1) by $(-q)^{r-c}$ and $(-q)^{c-r}$,
respectively:
$$\align X_{rc}A(r\,c) &\equiv - \sum_{j<c} (-q)^{j-c} X_{rj}
A(r\,j) \qquad\pmod{L} \tag2 \\
A(r\,c)X_{rc} &\equiv - \sum_{j<c} (-q)^{c-j}
A(r\,j) X_{rj} \qquad\pmod{L}. \tag3 \endalign$$

For any $j$, the third relation of \cite{\PaWa, Lemma 5.1.2} implies
that
$$X_{rj}A(r\,j) \equiv A(r\,j) X_{rj}+ (1-q^{-2}) \sum_{l<j} (-q)^{j-l}
A(r\,l) X_{rl} \qquad\pmod{L}, \tag4$$
since $X_{sj}A(s\,j) \in L$ for $s>r$.
Substituting (4) into (2) for all $j<c$, we obtain
$$\aligned X_{rc} A(r\,c) &\equiv - \sum_{j<c} (-q)^{j-c} A(r\,j)
X_{rj}\\
 &\qquad -(1-q^{-2}) \sum_{j<c} \sum_{l<j} (-q)^{2j-l-c} A(r\,l)
X_{rl} \qquad\pmod{L}\\
 &= - \sum_{l<c} \biggl[ (-q)^{l-c}+ (1-q^{-2}) \sum_{l<j<c}
(-q)^{2j-l-c} \biggr] A(r\,l) X_{rl}. \endaligned \tag5$$
The expression in square brackets can be simplified as follows:
$$\aligned (-q)^{l-c}+ (1-q^{-2}) \sum_{l<j<c} (-q)^{2j-l-c} &=
(-q)^{l-c} \biggl[ 1+ (1-q^{-2})\sum_{0<m<c-l} (-q)^{2m} \biggr]\\
 &= (-q)^{c-l-2}. \endaligned \tag6$$ 
Substituting (6) into (5) and replacing $l$ by $j$, we obtain
$$X_{rc} A(r\,c) \equiv - \sum_{j<c} (-q)^{c-j-2} A(r\,j) X_{rj}
\qquad\pmod{L}. \tag7$$ 
Finally, combining (3) with (7), we conclude that
$$[I\mid J] X_{rc}= A(r\,c) X_{rc} \equiv q^2 X_{rc} A(r\,c)= q^2
X_{rc} [I\mid J] \qquad\pmod{L},$$ 
as desired. \qed\enddemo

\proclaim{1.4\. Corollary} Let $I,J\subseteq \{1,\dots,n\}$ with
$|I|=|J|$, and set
$$L= \bigl\langle\; [I'\mid J'] \bigm| |I'|=|I|, \operatorname{\ and\ }
I'<I \operatorname{\ or\ } J'<J \;\bigr\rangle.$$
Then the coset $d= [I\mid J] +L$ generates a denominator set in $A/L$.
\endproclaim

\demo{Proof} Set $B= A/L$, and set $x_{ij}= X_{ij}+L$ for all $i,j$.
Lemma 1.3 says that
$$ dx_{ij}= q^{2-\delta(i,I)-\delta(j,J)} x_{ij}d \tag1$$
whenever $i\le \max(I)$ or $j\le \max(J)$. Hence, in this case
we have $d^rx_{ij} \in Bd^r$ for all $r\ge 0$.

When $i>\max(I)$ and $j>\max(J)$, Lemma 5.7 says that
$$dx_{ij}- q^2x_{ij}d= e := (1-q^2)[I\sqcup \{i\} \mid J\sqcup \{j\}]
+L. \tag2$$
Observe that $d$ and $e$ commute. Hence, it follows from (2) by an easy
induction that
$$\aligned d^rx_{ij} &= q^{2r}x_{ij}d^r+ (q^{2r-2}+ \cdots+ q^2+1)
ed^{r-1}\\
 &= \bigl[ q^{2r}x_{ij}d + (q^{2r-2}+ \cdots+ q^2+1)e \bigr] d^{r-1}
\endaligned \tag3$$
for all $r>0$. Combining (1) and (3), we see that
$$\xalignat2 d^rx_{ij} \ &\in\ Bd^{r-1} &&(i,j= 1,\dots,n;\
r=1,2,\dots). \tag4 \endxalignat$$

Since $B$ is spanned by products of the $x_{ij}$, it follows from (4)
that $D:= \{d^r \mid r\ge 0\}$ is a left Ore set in $B$. Similarly,
$D$ is right Ore, and therefore $D$ is a denominator set because $B$ is
noetherian. \qed\enddemo

\proclaim{1.5\. Proposition} Let $(\rc) \in \RC_t$, and set
$$L= \bigl\langle\; [I\mid J] \bigm| |I|=l\le t, \operatorname{\ and\ }
I< \{r_1,\dots,r_l\} \operatorname{\ or\ } J< \{c_1,\dots,c_l\}
\;\bigr\rangle.$$ 
Then the image of $\Drc$ in $A/L$ is a denominator set. Consequently,
$\Drctil$ is a denominator set in $A/\Krc$.
\endproclaim

\demo{Proof} For each $l=1,\dots,t$, Corollary 1.4 shows that
$\drc{l}+L$ generates a denominator set in $A/L$. The proposition
follows.
\qed\enddemo

\definition{1.6} In view of Proposition 1.5, we can form Ore
localizations
$$\Arc= (A/\Krc)[\Drctil^{-1}]$$
for $(\rc) \in \RC$. (It will follow from Lemma 2.5 that $\Arc \ne 0$.)
The localization maps $A\rightarrow A/\Krc
\rightarrow \Arc$ induce Zariski homeomorphisms
$$\specrc A \longrightarrow \spec \Arc.$$

\enddefinition

\proclaim{Lemma 1.7} Let $I,J \subseteq \{1,\dots,n\}$ with
$|I|=|J|$, and let $P$ be an ideal
of $A$.

{\rm (a)} Fix $J_1\subseteq J$. If $[I_1\mid J_1] \in P$ for all
$I_1\subseteq I$ with $|I_1|= |J_1|$, then $[I\mid J] \in P$.

{\rm (b)} Fix $I_1\subseteq I$. If $[I_1\mid J_1] \in P$ for all
$J_1\subseteq J$ with $|J_1|= |I_1|$, then $[I\mid J] \in P$.
\endproclaim

\demo{Proof} By symmetry (see (5.1)(b)), we need only prove (a).
Set $J_2= J\setminus J_1$. Then Lemma 5.4(a) provides a relation of the
form
$$\sum_{I_1\sqcup I_2=I} \pmqdot [I_1\mid J_1][I_2\mid J_2]= \pmqdot
[I\mid J],$$
where $\pm$ is an unspecified sign and $q^\bullet$ stands for an
unspecified power of
$q$. Since all the $[I_1\mid J_1] \in P$ by assumption, it follows that
$[I\mid J] \in P$. \qed\enddemo

\proclaim{Lemma 1.8} Let $I_1,I_2,J_1,J_2 \subseteq \{1,\dots,n\}$ with
$|I_1|=|J_1|$ and $|I_2|=|J_2|$, and let $P$ be a completely prime ideal
of $A$. Assume that one of the following conditions {\rm
(a)}, {\rm (b)}, or {\rm (c)} holds:
\roster
\item"(a)" {\rm (1)} $|I_1\cap I_2| >|J_1\cap J_2|$, and

\noindent {\rm (2)} $[I_1\mid J'] \in P$ whenever $J_1\cap J_2
\subseteq J'
\subseteq J_1\cup J_2$ with $|J'|=|J_1|$ but $J' \ne J_1$.
\endroster
\roster
\item"(b)" {\rm (1)} $|I_1\cap I_2| <|J_1\cap J_2|$, and

\noindent {\rm (2)} $[I'\mid J_1] \in P$ whenever $I_1\cap I_2
\subseteq I'
\subseteq I_1\cup I_2$ with $|I'|=|I_1|$ but $I' \ne I_1$.
\endroster
\roster
\item"(c)" {\rm (1)} $|I_1\cap I_2| =|J_1\cap J_2|$ and $[I_1\cup I_2
\mid J_1\cup J_2] \in P$, and

\noindent {\rm (2)} Either {\rm (a)(2)} or {\rm (b)(2)} holds.
\endroster
Then either $[I_1\mid J_1]\in P$ or $[I_2\mid J_2]\in P$.
\endproclaim

\demo{Proof} By symmetry, it suffices to prove cases (a) and (c).

(a) Set $V= J_1\cup J_2= J\sqcup L$ where $J= J_1\cap J_2$ and $L=
V\setminus J$. Since $|I_1\cap I_2| >|J|$, we have $|I_1\cup
I_2| <|V|$, and so there exists $U\subseteq \{1,\dots,n\}$ such that
$I_1\cup I_2 \subset U$ and $|U|=|V|$. Also, we have $|I_1|+|I_2|=
|J_1|+|J_2|= 2|J|+|L|$. Thus, Lemma 5.6(b) yields a relation of the form
$$\sum_{L= L'\sqcup L''} \pmqdot [I_1 \mid J \sqcup L'] [I_2 \mid
J\sqcup L'']=0. \tag\dagger$$
Now for each term in this sum, $J_1\cap J_2\subseteq J\sqcup L' \subseteq
J_1\cup J_2$ with $|J\sqcup L'|= |J_1|$. If $L' \ne J_1\setminus J$, then
$J\sqcup L' \ne J_1$, whence $[I_1\mid J\sqcup L'] \in P$ by hypothesis.
Therefore the remaining term in $(\dagger)$ must lie in $P$. This is the
term with $L'= J_1\setminus J$, whence $L''= J_2\setminus J$, and so
$J\sqcup L'= J_1$ and $J\sqcup L''= J_2$. Thus
$$\pmqdot [I_1\mid J_1] [I_2\mid J_2] \in P.$$
Since $P$ is completely prime, either $[I_1\mid J_1]\in P$ or $[I_2\mid
J_2]\in P$. 

(c) By symmetry, we may assume that (a)(2) holds. Again, set $V= J_1\cup
J_2= J\sqcup L$ where $J= J_1\cap J_2$ and $L= V\setminus J$. Set $U=
I_1\cup I_2$, and observe that  $|I_1|+|I_2|= 2|J|+|L|$. This time, Lemma
5.6(b) provides a
relation of the form
$$\sum_{L= L'\sqcup L''} \pmqdot [I_1 \mid J \sqcup L'] [I_2 \mid J\sqcup
L'']= \pmqdot [I_1\cap I_2\mid J] [U\mid V]. \tag\ddagger$$
Since $[U\mid V] \in P$ by hypothesis, the right hand side of
$(\ddagger)$ lies in $P$. Therefore we can proceed as in the proof
above. \qed\enddemo

\proclaim{1.9\. Theorem} Assume that $q\ne \pm1$. Let $P$ be a
completely prime ideal of $A$, and let $t\le n$ be maximal such
that $P$ does not contain all $t\times t$ quantum minors. Choose 
$$ \Itil = \{r_1< \cdots< r_t\} \subseteq \{1,\dots,n\}
\qquad\text{and}\qquad
\Jtil = \{c_1< \cdots< c_t\} \subseteq \{1,\dots,n\}$$
with $\Itil$ minimal such that some $[\Itil\mid *] \notin P$, and $\Jtil$
minimal such that some $[*\mid \Jtil] \notin P$. Then

{\rm (a)} $[r_1\cdots r_s\mid c_1\cdots c_s] \notin P$ for $s=1,\dots,t$.
In particular, $[\Itil \mid \Jtil] \notin P$.

{\rm (b)} $[I\mid J] \in P$ whenever $|I|=s\le t$ and either $I \ngeq
\{r_1,\dots,r_s\}$ or $J \ngeq \{c_1,\dots,c_s\}$. In particular, it
follows  that $\Itil$ and $\Jtil$ are unique.
\endproclaim

\demo{Proof} Since the theorem holds trivially when $t=0$, we may assume
that $t>0$.  
By assumption, there exists $J_0$ such that $[\Itil\mid J_0] \notin P$.
We first claim that
\roster
\item"(1)" $[I\mid J] \in P$ whenever $|I|=s\le t$ and $I<
\{r_1,\dots,r_s\}$.
\endroster
Suppose not, so that some $[I\mid J] \notin P$ where $|I|=s\le t$ and $I<
\{r_1,\dots,r_s\}$. We may assume that $s$ is minimal for this, and that
with $I$ fixed, $|J\cap J_0|$ is maximal. Note from the minimality of
$\Itil$ that $s<t$.

Write $I= \{i_1< \cdots< i_s\}$. There is some $b\le s$ such that
$i_l=r_l$ for $l<b$ while $i_b<r_b$, whence $\{i_1,\dots,i_b\} <
\{r_1,\dots,r_b\}$. Since $[I\mid J] \notin P$, Lemma 1.7(b) implies
that some $[i_1\cdots i_b \mid *] \notin P$. Hence, the minimality of
$s$ implies that $b=s$. Thus, we have $i_l=r_l$ for $l<s$ while
$i_s<r_s$. In particular, $r_{s-1}= i_{s-1}< i_s< r_s$.

Assume that $|J\cap J_0| \ge s-1$. In this case, we will apply Lemma
1.8 with
$$\xalignat 4 I_1 &= \Itil &I_2 &=I &J_1 &= J_0 &J_2 &= J. \endxalignat$$
Note that $|I_1\cap I_2|= s-1$ and $[I_1\mid J_1],\, [I_2\mid J_2] \notin
P$. When $I_1\cap I_2 \subseteq I' \subseteq I_1\cup I_2$ with $|I'|=t$
and $I'\ne I_1$, we have $I'= \Itil \cup \{i_s\} \setminus \{r_l\}$ for
some $l\ge s$. Then since $r_{s-1}< i_s< r_s \le r_l$, we have
$I'<\Itil$, and so
$[I'\mid J_1] \in P$ by the minimality of $\Itil$. Further, if $|J_1\cap
J_2|= s-1$, then $[I_1\cup I_2\mid J_1\cup J_2] \in P$ because $P$
contains all
$(t+1)\times (t+1)$ quantum minors. Therefore by Lemma 1.8(b), if
$|J\cap J_0| \ge s$, or by Lemma 1.8(c), if $|J\cap J_0| =s-1$, we have
either
$[I_1\mid J_1]\in P$ or
$[I_2\mid J_2]\in P$, giving us a contradiction. Therefore $|J\cap
J_0| <s-1$.

Next, we will apply Lemma 1.8 with the roles of $I_1,I_2$ and $J_1,J_2$
reversed, that is, with
$$\xalignat 4 I_1 &= I &I_2 &=\Itil &J_1 &= J &J_2 &= J_0. \endxalignat$$
When $J_1\cap J_2 \subseteq J' \subseteq J_1\cup J_2$ with $|J'|=s$ and
$J'\ne J_1$, we must have $|J'\cap J_0|= |J'\cap J_2| > |J\cap J_0|$. By
the maximality of $|J\cap J_0|$, we obtain $[I\mid J'] \in P$ in this
case. But then Lemma 1.8(a) leads to the same contradiction. 

Therefore (1)
holds. By symmetry, we must also have
\roster
\item"(2)" $[I\mid J] \in P$ whenever $|I|=s\le t$ and $J<
\{c_1,\dots,c_s\}$.
\endroster
We now proceed by induction on $s=1,\dots,t$ to verify the following
properties:
\roster
\item"($P_s$)" $[r_1\cdots r_s\mid c_1\cdots c_s] \notin P$;
\item"($Q_s$)" $[I\mid J]\in P$  whenever $|I|=s$ and either $I \ngeq
\{r_1,\dots,r_s\}$ or $J \ngeq \{c_1,\dots,c_s\}$.
\endroster
The theorem will then be established.

To start, note that we cannot have all $[r_1\mid *]\in P$ or all $[*\mid
c_1] \in P$, by Lemma 1.7. Choose
$i,j$ such that
$[r_1\mid j],\, [i\mid c_1]
\notin P$. If $j=c_1$ or $i=r_1$, then $[r_1\mid c_1] \notin P$.
Otherwise, in view of (1) and (2) we must have $j>c_1$ and $i>r_1$.
Hence, because of the assumption that $q\ne \pm1$, we have
$$[r_1\mid c_1] [i\mid j]- [i\mid j] [r_1\mid c_1] = (q-q^{-1}) [r_1\mid
j] [i\mid c_1] \notin P,$$
which implies that $[r_1\mid c_1] \notin P$. Therefore ($P_1$) holds.
Property ($Q_1$) is immediate from (1) and (2).

Now let $1<s\le t$ and assume that ($P_a$) and ($Q_a$) hold for all
$a<s$. By Lemma 1.7, there exist $j_1< \cdots< j_s$ such that
$$[r_1\cdots r_s\mid j_1\cdots j_s] \notin P,$$
and we may assume that
$\{j_1,\dots,j_s\}$ is minimal for this. Likewise, there exist $i_1<
\cdots< i_s$ such that
$[i_1\cdots i_s\mid c_1\cdots c_s] \notin P$ and such that
$\{i_1,\dots,i_s\}$ is minimal for this. We cannot have
$\{i_1,\dots,i_{s-1}\} \ngeq \{r_1,\dots,r_{s-1}\}$, since then
($Q_{s-1}$) would imply that all $[i_1\cdots i_{s-1} \mid *] \in P$,
whence Lemma 1.7(b) would imply that all $[i_1\cdots i_s\mid *] \in P$.
Therefore $\{i_1,\dots,i_{s-1}\} \ge \{r_1,\dots,r_{s-1}\}$, and similarly
$\{j_1,\dots,j_{s-1}\} \ge \{c_1,\dots,c_{s-1}\}$.

Suppose there exists $b<s$ such that $j_l=c_l$ for $l<b$ while $j_b>c_b$.
We will apply Lemma 1.8 with
$$\xalignat2 I_1 &= \{r_1,\dots,r_s\} &I_2 &= \{r_1,\dots,r_b\} \\
J_1 &= \{j_1,\dots,j_s\} &J_2 &= \{c_1,\dots,c_b\}. \endxalignat$$
Observe that
$|I_1\cap I_2| =b> |J_1\cap J_2|$. If $J_1\cap J_2 \subseteq J' \subseteq
J_1\cup J_2$ with $|J'|=s$ and
$J'\ne J_1$, then $J'= J_1\cup \{c_b\} \setminus \{j_d\}$ for some $d\ge
b$. In this case, $J'< \{j_1,\dots,j_s\}$, whence $[I_1\mid J'] \in P$
by the minimality of $\{j_1,\dots,j_s\}$. Therefore Lemma 1.8(a)
implies that either $[I_1\mid J_1]\in P$ or $[I_2\mid J_2]\in P$. But
$[I_1\mid J_1] \notin P$ by choice of $J_1$, and $[I_2\mid J_2] \notin
P$ by ($P_b$), so we have a contradiction. Therefore $j_l=c_l$ for all
$l<s$. Similarly, $i_l=r_l$ for all $l<s$.

If $j_s=c_s$, then ($P_s$) holds. If $j_s<c_s$, then $\{j_1,\dots,j_s\}< 
\{c_1,\dots,c_s\}$, which would imply $[r_1\cdots r_s\mid j_1\cdots j_s]
\in P$ by (2), contradicting our assumptions. Therefore we may assume
that $j_s>c_s$. Likewise, we may assume that $i_s>r_s$.

Set $U= \{r_1,\dots,r_s,i_s\}$ and $V= 
\{c_1,\dots,c_s,j_s\}$. By Lemma 5.3(d), we have
$$\multline [i_1\cdots i_s\mid j_1\cdots j_s] [r_1\cdots r_s\mid
c_1\cdots c_s] - [r_1\cdots r_s\mid c_1\cdots c_s] [i_1\cdots i_s\mid
j_1\cdots j_s] \\
 = (q^{-1}-q) [r_1\cdots r_s\mid j_1\cdots j_s] 
[i_1\cdots i_s\mid c_1\cdots c_s]. \endmultline$$
Since neither of the factors $[r_1\cdots r_s\mid j_1\cdots j_s]$ and
$[i_1\cdots i_s\mid c_1\cdots c_s]$ is in
$P$, it follows that
$[r_1\cdots r_s\mid c_1\cdots c_s]$ cannot be in $P$. This establishes
property ($P_s$).

Finally, suppose that ($Q_s$) fails. By symmetry, we may assume that
$[I\mid J] \notin P$ for some $I,J$ with $|I|=s$ and $I \ngeq
\{r_1,\dots,r_s\}$. We may also assume that $I$ is minimal for this, and
that with $I$ fixed, $J$ is minimal.

Write $I= \{i_1< \cdots< i_s\}$ and $J= \{j_1< \cdots< j_s\}$. If
$\{i_1,\dots,i_{s-1}\} \ngeq \{r_1,\dots,r_{s-1}\}$, then by ($Q_{s-1}$)
we would have all $[i_1\cdots i_{s-1} \mid *] \in P$, whence Lemma
1.7(b) would imply that all $[I\mid *] \in P$, contradicting our choice
of $I$. Thus $\{i_1,\dots,i_{s-1}\} \ge \{r_1,\dots,r_{s-1}\}$, and
similarly
$\{j_1,\dots,j_{s-1}\} \ge \{c_1,\dots,c_{s-1}\}$. Since $I\ngeq
\{r_1,\dots,r_s\}$, we must also have $i_s<r_s$. Note that $r_{s-1}\le
i_{s-1}< i_s< r_s$, and so $i_s\notin \{r_1,\dots,r_{s-1}\}$. Further,
$\{r_1,\dots,r_{s-1},i_s\}< \{r_1,\dots,r_s\}$, and so all $[r_1\cdots
r_{s-1}i_s \mid *] \in P$ by (1), whence $I\ne
\{r_1,\dots,r_{s-1},i_s\}$. Therefore there is some $b<s$ such that
$i_l=r_l$ for $l<b$ while $i_b>r_b$.

Suppose there exists $d\le b$ such that $j_m=c_m$ for $m<d$ while
$j_d>c_d$. We will apply Lemma 1.8 with
$$\xalignat2 I_1 &= I= \{i_1,\dots,i_s\} &I_2 &= \{r_1,\dots,r_d\}\\
J_1 &= J= \{j_1,\dots,j_s\} &J_2 &= \{c_1,\dots,c_d\}. \endxalignat$$
Note that $|I_1\cap I_2| \ge d-1= |J_1\cap J_2|$, and that $|I_1\cap I_2|
= |J_1\cap J_2|$ only when $d=b$. Since $r_{b-1}< r_b< i_b$, we have 
$\{r_1,\dots,r_b,i_b,\dots,i_{s-1}\} <I$, and so all $[r_1\cdots
r_bi_b\cdots i_{s-1} \mid *] \in P$ by the minimality of $I$. Then
Lemma 1.7(b) implies that all $[r_1\cdots
r_bi_b\cdots i_s \mid *] \in P$. In particular, when $d=b$ we find that
$[I_1\cup I_2\mid J_1\cup J_2] \in P$.

If $J_1\cap J_2 \subseteq J' \subseteq
J_1\cup J_2$ with $|J'|=s$ and
$J'\ne J_1$, then $J'= J\cup \{c_d\} \setminus \{j_p\}$ for some $p\ge
d$. In this case, $J'<J$, and so $[I_1\mid J'] \in P$ by the minimality
of $J$. Hence, case (a) of Lemma 1.8 (if $d>b$) or case (c) (if $d=b$)
implies that either
$[I_1\mid J_1]\in P$ or
$[I_2\mid J_2]\in P$. But $[I_1\mid J_1] \notin P$ by assumption,
and $[I_2\mid J_2] \notin P$ by ($P_d$), so we have a contradiction.
Therefore $j_m=c_m$ for all $m\le b$.

We will conclude by applying Lemma 1.8 with
$$\xalignat2 I_1 &= I= \{i_1,\dots,i_s\} &I_2 &= \{r_1,\dots,r_b\}\\ 
J_1 &= J= \{j_1,\dots, j_s\} &J_2 &= \{c_1,\dots,c_b\}. \endxalignat$$
Note that $i_{b-1}=r_{b-1}< r_b< i_b$ implies $r_b \notin I$, and so
$|I_1\cap I_2| <|J_1\cap J_2|$. If $I_1\cap I_2 \subseteq I' \subseteq
I_1\cup I_2$ with $|I'|=s$ and $I'\ne I_1$, then $I'= I\cup \{r_b\}
\setminus \{i_p\}$ for some $p\ge b$. In this case, $I'<I$, whence
$[I'\mid J_1] \in P$ by the minimality of $I$. Thus Lemma 1.8(b)
implies  that either $[I_1\mid J_1]\in P$ or
$[I_2\mid J_2]\in P$, and again we have reached a contradiction.

Therefore ($Q_s)$ must hold, which establishes our induction step.
\qed\enddemo

\proclaim{1.10\. Corollary} Assume that $q\ne \pm1$.
Given any completely prime ideal $P\in \spec A$, there is a unique pair
$(\rc)
\in \RC$ such that $\Krc \subseteq P$ and $P\cap \Drc = \varnothing$.

Thus, if $q$ is not a root of unity,
$$\spec A = \bigsqcup_{(\rc) \in \RC} \specrc A.$$
\endproclaim

\demo{Proof} Let $t\le n$ be maximal such
that $P$ does not contain all $t\times t$ quantum minors, let
$\{r_1< \cdots< r_t\}$ and $\{c_1< \cdots< c_t\}$ be as in Theorem
1.9, and set $\bfr= (r_1,\dots,r_t)$ and $\bfc= (c_1,\dots,c_t)$. The
theorem implies that $\Krc
\subseteq P$ and that $\drc{s}
\notin P$ for $s=1,\dots,t$. Since $P$ is completely prime, it follows
that $P\cap \Drc = \varnothing$.

Now suppose that we also have $(\bfr',\bfc')
\in \RC_{t'}$ for some $t'$ such that $K_{\bfr',\bfc'} \subseteq P$ and
$P\cap D_{\bfr',\bfc'} =
\varnothing$. Then $P$ contains all $(t'+1) \times (t'+1)$ quantum
minors but not all $t'\times t'$ quantum minors, whence $t'=t$.
Moreover, we have $\drc{t} \notin K_{\bfr',\bfc'}$ and
$d_t^{\bfr',\bfc'}
\notin \Krc$. The first relation implies that $\bfr \ge \bfr'$ and
$\bfc \ge \bfc'$, and the second relation yields the reverse
inequalities. (Here we have transferred the relation $\le$ in (1.1)
from index sets to sequences in the obvious manner.) Therefore
$\bfr'= \bfr$ and
$\bfc'= \bfc$.
\qed\enddemo

\head 2.  Structure of $\Arc$ \endhead

The purpose of this section is to develop a structure theorem for the
localizations $\Arc$. We introduce localized factor algebras $\Rrplus$
and $\Rcminus$ of $A$ (patterned after quantized coordinate rings of
groups of triangular matrices) together with subalgebras $\Brplus
\subset \Rrplus$ and $\Bcminus \subset \Rcminus$ (patterned after
quantized coordinate rings of unipotent groups of triangular matrices),
and we show that $\Arc$ is isomorphic to an algebra $\Brc$ trapped
between $\Brplus\otimes\Bcminus$ and $\ror$. More precisely, we prove
that $\Brc$ is a skew-Laurent extension of $\Brplus\otimes\Bcminus$
(this is a key ingredient in establishing that $\Arc\cong\Brc$), and
that $\ror$ is a skew-Laurent extension of $\Brc$.

\definition{2.1} Fix $t\in \{0,1,\dots,n\}$ and $(\rc) \in \RC_t$
throughout the section. Set 
$$\Rrpluszero= A/\langle X_{ij}\mid j>t \text{\ or\ } i<r_j \rangle
\qquad\text{and}\qquad \Rcminuszero= A/\langle X_{ij}\mid i>t \text{\ or\ }
j<c_i \rangle.$$
 Write
$Y_{ij}$ and $Z_{ij}$ for the images of $X_{ij}$ in
$\Rrpluszero$ and $\Rcminuszero$, respectively. Note that these
algebras are iterated skew polynomial extensions of $k$, hence
noetherian domains, the natural indeterminates for these iterated skew
polynomial structures being those $Y_{ij}$ and $Z_{ij}$ which are
nonzero. These indeterminates, when recorded within an $n\times n$
matrix, display `stairstep' patterns -- for example, if $n=4$ and
$\bfr= (1,2,4)$, the $Y_{ij}$ may be displayed as follows:
$$\left[ \smallmatrix Y_{11} &0&0&0\\ Y_{21} &Y_{22} &0&0\\ Y_{31}
&Y_{32} &0&0\\ Y_{41} &Y_{42} &Y_{43} &0\endsmallmatrix \right].$$
Since all the information is recorded in the placement of the zero and
nonzero positions within this matrix, a convenient abbreviation for
this example is to write
$$R^+_{(1,2,4),0}= k \left[ \smallmatrix
\scriptscriptstyle{+} &\scriptscriptstyle{0} &\scriptscriptstyle{0}
&\scriptscriptstyle{0}\\
\scriptscriptstyle{+} &\scriptscriptstyle{+} &\scriptscriptstyle{0}
&\scriptscriptstyle{0}\\
\scriptscriptstyle{+} &\scriptscriptstyle{+} &\scriptscriptstyle{0}
&\scriptscriptstyle{0}\\
\scriptscriptstyle{+} &\scriptscriptstyle{+} &\scriptscriptstyle{+}
&\scriptscriptstyle{0} \endsmallmatrix   \right] .$$

Observe that the
$Y_{r_ss}$ are regular normal elements
in $\Rrpluszero$, and that the $Z_{sc_s}$
are regular normal elements in $\Rcminuszero$. More precisely,
$$\align Y_{r_ss}Y_{ij} &= 
\cases q^{-1}Y_{ij}Y_{r_ss} &\qquad (i=r_s,\, j\ne s)\\ 
qY_{ij}Y_{r_ss} &\qquad (i\ne r_s,\, j=s)\\ 
Y_{ij}Y_{r_ss} &\qquad (i\ne r_s,\, j\ne s) \endcases \\ 
Z_{sc_s}Z_{lm} &=
\cases qZ_{lm}Z_{sc_s} &\qquad (l=s,\, m\ne c_s)\\
q^{-1}Z_{lm}Z_{sc_s} &\qquad (l\ne s,\, m=c_s)\\
Z_{lm}Z_{sc_s} &\qquad (l\ne s,\, m\ne c_s). \endcases \endalign$$
(For instance, the first relation above holds when $j>s$ because
$Y_{ij}=0$ in that case. To verify the third relation, observe that
$Y_{r_sj}=0$ if $j>s$, while $Y_{is}=0$ if $i<r_s$.) In particular,
the $Y_{r_ss}$ commute with each other, and the
$Z_{sc_s}$ commute with each other.

Due to the normality of the $Y_{r_ss}$ and the $Z_{sc_s}$, we can form
Ore localizations
$$\Rrplus= \Rrpluszero[ Y^{-1}_{r_11}, Y^{-1}_{r_22}, \dots,
Y^{-1}_{r_tt}] \qquad\text{and}\qquad \Rcminus= \Rcminuszero[
Z^{-1}_{1c_1}, Z^{-1}_{2c_2}, \dots, Z^{-1}_{tc_t}].$$
These algebras are noetherian domains, and they may be viewed as
quantized coordinate rings of certain locally closed subvarieties of
$M_n(k)$. Extending the abbreviated description given for the example
above, we display the following abbreviation for
$\Rrplus$ in that case:
$$R^+_{(1,2,4)}= k \left[ \smallmatrix
\scriptscriptstyle{\pm} &\scriptscriptstyle{0} &\scriptscriptstyle{0}
&\scriptscriptstyle{0}\\
\scriptscriptstyle{+} &\scriptscriptstyle{\pm} &\scriptscriptstyle{0}
&\scriptscriptstyle{0}\\
\scriptscriptstyle{+} &\scriptscriptstyle{+} &\scriptscriptstyle{0}
&\scriptscriptstyle{0}\\
\scriptscriptstyle{+} &\scriptscriptstyle{+} &\scriptscriptstyle{\pm}
&\scriptscriptstyle{0} \endsmallmatrix   \right] .$$
Observe that the standard $\ZZ^{2n}$-grading on $A$ induces
$\ZZ^{2n}$-gradings on $\Rrplus$ and $\Rcminus$, which we also refer
to as {\it standard\/}.
\enddefinition

\definition{2.2} Let $\pirpluszero: A\rightarrow \Rrpluszero$ and
$\picminuszero: A\rightarrow \Rcminuszero$ be the quotient maps, and
define
$$\betarc: A @>{\Delta}>> A\otimes A
@>{\pirpluszero\otimes\picminuszero}>> \Rrpluszero\otimes \Rcminuszero
@>{\subseteq}>> \ror.$$
Observe that
$$\betarc(X_{ij})= \sum_{l\le t,\, r_l\le i,\, c_l\le j} Y_{il}\otimes
Z_{lj}$$
for all $i,j$. In particular, $\betarc(X_{ij})=0$ when $i<r_1$ or
$j<c_1$.
\enddefinition

\proclaim{2.3\. Lemma} $\Krc\subseteq \ker(\betarc)$. \endproclaim

\demo{Remark} We conjecture that $\ker(\betarc)= \Krc$.
\enddemo

\demo{Proof} Since $X_{ij}\in \ker(\pirpluszero)$ for  $j>t$, all
$(t+1)\times (t+1)$ and larger quantum minors lie in
$\ker(\pirpluszero)$. In view of the rule for comultiplication of
quantum minors (see (5.1)(d)), it follows that all $(t+1)\times (t+1)$
and larger quantum minors lie in $\ker(\betarc)$.

Consider an index set $I$ with $|I|=l\le t$ and $I \not\ge
\{r_1,\dots,r_l\}$. Write $I= \{i_1< \cdots< i_l\}$; then $i_m<r_m$
for some $m\le l$. Hence, $Y_{i_sj}=0$ for all $s\le m$ and $j\ge
m$. This implies that $\pirpluszero[i_1\cdots i_m\mid M] =0$ for all
$M$ with
$|M|=m$, whence $\pirpluszero[I\mid K]=0$ for all $K$ with $|K|=l$
(cf\. Lemmas 1.7 or 5.4). Therefore
$\betarc[I\mid J]=0$ for all $J$ with $|J|=l$.

Likewise, $\betarc[I\mid J]=0$ whenever $|I|=l\le t$ and $J \not\ge
\{c_1,\dots,c_l\}$. \qed\enddemo

\definition{2.4} Let $\Brc$ denote the $k$-subalgebra of $\ror$
generated by the set
$$\{Y_{il}\otimes Z_{lj} \mid l\le t,\ i\ge r_l,\ j\ge c_l\} \cup
\{Y_{r_ll}^{-1} \otimes Z_{lc_l}^{-1} \mid l\le t\}.$$
We may also express $\Brc$ as the subalgebra of $\ror$ generated by 
$$\{ Y_{il}Y_{r_ll}^{-1} \otimes 1 \mid l\le t,\ i> r_l \} \cup
\{ 1\otimes Z_{lj}Z_{lc_l}^{-1} \mid l\le t,\ j> c_l \} \cup
\{ (Y_{r_ll} \otimes Z_{lc_l})^{\pm1} \mid l\le t\}.$$
Note that $\betarc(X_{ij}) \in \Brc$ for all $i,j$, so that
$\betarc(A) \subseteq \Brc$.

Let $l\le t$. Since
$Y_{r_sj}=0$ for
$s\le l$ and
$j>s$, we have
$$\pirpluszero[r_1 \cdots r_l\mid K]= \cases Y_{r_11} Y_{r_22}
\cdots Y_{r_ll} &\quad (K= \{1,\dots, l\}) \\
0 &\quad (K\ne \{1,\dots, l\}). \endcases$$
Similarly,
$\picminuszero[1\cdots l\mid c_1\cdots c_l]= Z_{1c_1} Z_{2c_2} \cdots
Z_{lc_l}$, and therefore
$$\betarc(\drc{l})= (Y_{r_11}\otimes Z_{1c_1}) (Y_{r_22}\otimes
Z_{2c_2})
\cdots (Y_{r_ll}\otimes Z_{lc_l}).$$
In particular, $\betarc(\drc{l})$ is invertible in $\Brc$. 
\enddefinition

\proclaim{2.5\. Lemma} The map $\betarc$ induces a surjective
$k$-algebra homomorphism
$$\betatilrc : \Arc \longrightarrow \Brc.$$
\endproclaim

\demo{Remark} We shall prove later (Theorem 2.11) that $\betatilrc$ is
an isomorphism. Note that the lemma already implies that $\Arc$ is
nonzero.
\enddemo

\demo{Proof} We have $\Krc \subseteq \ker(\betarc)$ by Lemma 2.3 and
$\betarc(A) \subseteq \Brc$ by (2.4), and so $\betarc$ induces a
homomorphism
$\betarc' : A/\Krc \rightarrow \Brc$. It also follows from (2.4) that
$\betarc'$ sends the elements of $\Drctil$ to units of $\Brc$, and
therefore $\betarc'$ does induce a homomorphism $\betatilrc : \Arc
\rightarrow \Brc$.

Set $E= \betatilrc(\Arc)$; we must show that $E= \Brc$. Note that
$$(Y_{r_ll} \otimes Z_{lc_l})^{\pm1} = \betarc(\drc{l})^{\pm1}
\betarc(\drc{l-1})^{\mp1} \in E$$
for $l\le t$, where $\drc{0}=1$. It remains to show
that $Y_{il}\otimes Z_{lj} \in E$ for all $i,l,j$.

As in (2.4), we see that
$$\align \betarc \bigl( [r_1\cdots r_{l-1} i\mid c_1\cdots c_l] \bigr)
&= (Y_{r_11}\otimes Z_{1c_1})\cdots (Y_{r_{l-1},l-1}\otimes
Z_{l-1,c_{l-1}}) (Y_{il}\otimes Z_{lc_l})\\
 &= \betarc(\drc{l-1}) (Y_{il}\otimes Z_{lc_l}) \endalign$$ 
for $l\le t$ and $i\ge r_l$, whence $Y_{il}\otimes Z_{lc_l} \in E$.
Similarly,
$Y_{r_ll}\otimes Z_{lj} \in E$ for $j\ge c_l$, and
therefore
$$Y_{il}\otimes Z_{lj}= q^{1-\delta(j,c_l)} (Y_{r_ll}^{-1} \otimes
Z_{lc_l}^{-1}) (Y_{r_ll}\otimes Z_{lj}) (Y_{il}\otimes Z_{lc_l}) \in
E,$$ as desired. \qed\enddemo

\proclaim{2.6\. Lemma} $\ror$ is a skew-Laurent extension of $\Brc$ of
the form
$$\ror= \Brc[ 1\otimes Z_{1c_1}^{\pm1},\, \dots,\, 1\otimes
Z_{tc_t}^{\pm1};\, \tau_1,\dots,\tau_t]$$
for some $\tau_1,\dots,\tau_t \in \{1\}^{2n}\times H$. \endproclaim

\demo{Proof} First observe that there exist $\sigma_1,\dots,\sigma_t
\in H$ such that $Z_{lc_l}r= \sigma_l(r)Z_{lc_l}$ for $l\le t$ and
$r\in \Rcminuszero$. This relation extends to $r\in \Rcminus$, and so
$$(1\otimes Z_{lc_l})w= (1\times\sigma_l)(w) (1\otimes Z_{lc_l})$$
for $l\le t$ and $w\in \ror$. In particular, if $\tau_l=
(1,\dots,1,\sigma_l)$ then $1\otimes Z_{lc_l}$ is $\tau_l$-normal with
respect to $\Brc$.

The standard  $\ZZ^{2n}$-gradings on $\Rrplus$ and $\Rcminus$
(cf\. (2.1)) induce a
$\ZZ^{4n}$-grading on
$\ror$. With respect to this grading, $\Brc$ is a
homogeneous subalgebra of $\ror$, and its homogeneous components
have degrees of the form $(*,b,b,*)$. On the other hand, $1\otimes
Z_{lc_l}$ has degree $(0,0,\epsilon_l,*)$, so the monomials
$$(1\otimes Z_{1c_1})^{m_1} (1\otimes Z_{2c_2})^{m_2} \cdots
(1\otimes Z_{tc_t})^{m_t}$$
have degrees $(0,0,m_1\epsilon_1 +\dots+ m_t\epsilon_t,*)$. It
follows that these monomials are left (or right) linearly
independent over $\Brc$. Hence, the subalgebra
$$C= \sum_{m_1,\dots,m_t\in \ZZ} \Brc(1\otimes Z_{1c_1})^{m_1}
(1\otimes Z_{2c_2})^{m_2} \cdots (1\otimes Z_{tc_t})^{m_t}$$ 
of $\ror$ is a skew-Laurent extension of $\Brc$ of the desired form.

It remains to show that $C=\ror$. First note that $Y_{il}\otimes 1=
(Y_{il}\otimes Z_{lc_l}) (1\otimes Z_{lc_l}^{-1}) \in C$ for $l\le t$
and $i\ge r_l$, and that
$$Y_{r_ll}^{-1}\otimes1= (Y_{r_ll}^{-1}\otimes Z_{lc_l}^{-1})
(1\otimes Z_{lc_l}) \in C$$
for $l\le t$. On the other hand,
$$1\otimes Z_{lj}= (Y_{r_ll}\otimes Z_{lj}) (Y_{r_ll}^{-1} \otimes
Z_{lc_l}^{-1}) (1\otimes Z_{lc_l}) \in C$$
for $l\le t$ and $j\ge c_l$, and $1\otimes Z_{lc_l}^{-1} \in C$ for
$l\le t$. Therefore $C=\ror$. \qed\enddemo

\definition{2.7} Let $\Brplus$ and $\Bcminus$ denote the subalgebras
of $\Rrplus$ and $\Rcminus$ generated by the respective subsets
$$\{ Y_{il}Y_{r_ll}^{-1} \mid l\le t,\ i> r_l \} \qquad\text{and}\qquad
\{ Z_{lj}Z_{lc_l}^{-1} \mid l\le t,\ j> c_l \}.$$
The algebra $\Brplus$, for instance, may be viewed as a quantized
coordinate ring of the variety
$$\{ (a_{ij})\in M_n(k) \mid a_{ij}=0 \text{\ when\ } j>t \text{\ or\ }
i<r_j, \text{\ and\ } a_{r_jj}=1 \text{\ for\ } j\le t\}.$$
(The factor algebra $\Rrplus/ \langle Y_{r_11}-1,\, \dots,\, Y_{r_tt}-1
\rangle$, which one might expect to appear in the above role, is
inappropriate because it collapses to
$k[Y_{r_11}^{\pm1},\dots,Y_{r_tt}^{\pm1}]$ when $q\ne1$.)

As noted in (2.4), $\Brc$ is generated by its subalgebra $\Brplus \otimes
\Bcminus$ together with the elements $(Y_{r_ll} \otimes
Z_{lc_l})^{\pm1}$ for
$l\le t$. In fact:
\enddefinition

\proclaim{Lemma} $\Brc$ is a skew-Laurent extension of $\Brplus
\otimes \Bcminus$ of the form
$$\Brc= (\Brplus \otimes \Bcminus) [(Y_{r_11} \otimes
Z_{1c_1})^{\pm1}, \dots, (Y_{r_tt} \otimes Z_{tc_t})^{\pm1};
\eta_1,\dots,\eta_t]$$
for some $\eta_1,\dots,\eta_t \in H\times H$. \endproclaim

\demo{Proof} This is proved in the same manner as Lemma 2.6.
\qed\enddemo

\definition{2.8} We use the above structure of $\Brc$ in constructing
a homomorphism
$\Brc \rightarrow \Arc$ which will be the inverse of $\betatilrc$.
To begin the construction, we will define suitable homomorphisms
from $\Brplus$ and $\Bcminus$ to $\Arc$ whose images centralize each
other. For that purpose, we need to know the defining relations for
$\Brplus$ and $\Bcminus$, as well as  certain commutation relations in
$\Arc$.

Set $y_{ij}= Y_{ij}Y_{r_jj}^{-1}$ for $j\le t$ and
$i>r_j$. In view of the basic commutation relations satisfied by the
$Y_{lm}$, it is easily checked that the $y_{ij}$ satisfy the following
relations:
$$\alignedat2 y_{ij}y_{im} &= qy_{im}y_{ij} &&(j<m)\\
y_{ij}y_{lj} &= qy_{lj}y_{ij} &&(i<l)\\
y_{ij}y_{lm} &= y_{lm}y_{ij} &&(i<l,\ j>m)\\
y_{ij}y_{lm} &= \cases y_{lm}y_{ij} &\quad (i<r_m)\\
q^{-1}y_{lm}y_{ij} +\qhat y_{lj} &\quad (i=r_m)\\
y_{lm}y_{ij}+ \qhat y_{im}y_{lj} &\quad (i>r_m) \endcases
 &\qquad\qquad&(i<l,\ j<m). \endalignedat \tag1$$
\enddefinition

\proclaim{Lemma} The relations {\rm (1)} are defining relations for
the elements $y_{ij}$ generating the algebra $\Brplus$. \endproclaim

\demo{Proof} Let $S$ be the $k$-algebra presented by generators
$s_{ij}$  for $j\le t$ and
$i>r_j$ satisfying the analogs of (1). Then there is a $k$-algebra
homomorphism $\phi: S\rightarrow \Brplus$ such that $\phi(s_{ij})=
y_{ij}$ for all $i,j$. 

In $\Rrplus$, list the $Y_{ij}$ lexicographically, and observe that
the ordered monomials in the $Y_{ij}$ are linearly independent. This
remains true for ordered monomials in which we allow the $Y_{r_jj}$ to
have negative exponents. Since the $Y_{r_jj}$ commute up to scalars
with the $Y_{lm}$ (recall (2.1)), it follows that the ordered
monomials in the $y_{ij}$ are linearly independent, and so these
monomials form a basis for $\Brplus$. On the other hand, there are
sufficient commutation relations for the $s_{ij}$ to show that the
ordered monomials in the $s_{ij}$ span $S$. Hence, $\phi$ maps a
spanning subset of $S$ to a basis for $\Brplus$. Therefore $\phi$ is
an isomorphism, and the lemma is proved. \qed\enddemo

\proclaim{2.9\. Lemma} Let $M= [I\sqcup \{a\} \mid J]$ and $N= [I'
\mid J']$ be quantum minors in $A$ with $I\subset I'$ and $J\subseteq
J'$. Assume that $a\notin I'$ and that $b= \max(I') \notin I$.

{\rm (a)} If $a>b$, then $MN= q^{-1}NM$.

{\rm (b)} If $a<b$, then 
$$MN- q^{-1}NM \equiv \qhat (-q)^{
|(I'\setminus I) \cap (a,b)| } [I\sqcup \{b\} \mid J] [I'\sqcup \{a\}
\setminus \{b\} \mid J']$$
modulo the ideal $L := \bigl\langle\, [I'\sqcup \{a\} \setminus \{i'\}
\mid J'] \bigm| i'\in I'\cap(a,b) \,\bigr\rangle$.
\endproclaim

\demo{Proof} (a) Expand $M$ using the $q$-Laplace relation of
Corollary 5.5(b2) with $r=a$. Since $(a,n] \cap I= \varnothing$, we
get
$$M= \sum_{j\in J} (-q)^{ |(j,n]\cap J| } [I \mid J\setminus \{j\}]
X_{aj}. \tag1$$
Note that all the $[I \mid J\setminus \{j\}]$ commute with $N$. Lemma
5.2(c1) implies that $X_{aj}N= q^{-1}NX_{aj}$ for all $j\in J$, and
thus it follows from (1) that $MN= q^{-1}NM$.

(b) Set $\alpha= |(a,n] \cap I|$, and again expand $M$ using Corollary
5.5(b2) with $r=a$. This time, we get
$$(-q)^\alpha M= \sum_{j\in J} (-q)^{ |(j,n]\cap J| } [I \mid
J\setminus \{j\}] X_{aj}. \tag2$$
For all $j\in J$, Lemma 5.2(c2) implies that
$$NX_{aj}- qX_{aj}N= \qhat \sum \Sb i'\in I'\\ i'>a \endSb (-q)^{
|I'\cap [a,i']| } X_{i'j} [I'\sqcup \{a\} \setminus \{i'\} \mid J'],$$
whence
$$X_{aj}N - q^{-1}NX_{aj} \equiv \qhat (-q)^\beta X_{bj} [I'\sqcup
\{a\} \setminus \{b\} \mid J'] \qquad \pmod{L}, \tag3$$
where $\beta= |I'\cap [a,b]| -1= |I' \cap (a,b)|$. Note that $\beta-
\alpha= |(I'\setminus I) \cap (a,b)|$.

Combining (2) and (3), we obtain
$$MN \equiv \sum_{j\in J} (-q)^{ |(j,n]\cap J| -\alpha } [I \mid
J\setminus \{j\}] \biggl( q^{-1}NX_{aj}+ \qhat (-q)^\beta X_{bj}
[I'\sqcup \{a\} \setminus \{b\} \mid J'] \biggr) \tag4$$
modulo $L$.
Since all the $[I \mid J\setminus \{j\}]$ commute with $N$, it follows
from (2) that
$$\sum_{j\in J} (-q)^{ |(j,n]\cap J| -\alpha } [I \mid
J\setminus \{j\}] q^{-1}NX_{aj}= q^{-1}NM. \tag5$$
Further, an application of Corollary 5.5(b2) with $r=b$ yields
$$\sum_{j\in J} (-q)^{ |(j,n]\cap J| } [I \mid
J\setminus \{j\}] X_{bj}= [I\sqcup \{b\} \mid J]. \tag6$$
Combining (4), (5), and (6), we complete the proof. \qed\enddemo

\proclaim{2.10\. Corollary} Let $M= [I\mid J\sqcup \{a\}]$ and $N=
[I'\mid J']$ be quantum minors in $A$ with $I\subseteq I'$ and
$J\subset J'$. Assume that $a\notin J'$ and that $b= \max(J') \notin
J$.

{\rm (a)} If $a>b$, then $MN= q^{-1}NM$.

{\rm (b)} If $a<b$, then 
$$MN- q^{-1}NM \equiv \qhat (-q)^{
|(J'\setminus J) \cap (a,b)| } [I \mid J\sqcup \{b\}] [I' \mid
J'\sqcup \{a\} \setminus \{b\}]$$
modulo the ideal $L := \bigl\langle\, [I' \mid J'\sqcup \{a\} \setminus
\{j'\}] \bigm| j'\in J'\cap(a,b) \,\bigr\rangle$.
\endproclaim

\demo{Proof} This follows from Lemma 2.9 by symmetry. \qed\enddemo

\proclaim{2.11\. Theorem} The map $\betatilrc : \Arc \longrightarrow
\Brc$ is an isomorphism. \endproclaim

\demo{Proof} Recall the notation $\drctil{l}= \drc{l}+\Krc$ from
(1.2). Similarly, we shall use tildes to denote other cosets in
$A/\Krc$. To abbreviate the relation of
congruence modulo $\Krc$, we adopt the notation $\congmodKrc$.
We shall use the same symbols for elements of $A/\Krc$ and their
images in the localization $\Arc$, which does not cause problems as
long as we only transfer equations from $A/\Krc$ to $\Arc$ and not in
the reverse direction.

We proceed to construct, in several steps, a $k$-algebra homomorphism
$\phi: \Brc \rightarrow
\Arc$ that will be an inverse for $\betarc$. The construction of $\phi$
is based on the skew-Laurent structure of
$\Brc$ given in Lemma 2.7. The first ingredient will be a homomorphism
from
$\Brplus$ to $\Arc$. To describe it, set
$$\xalignat2 u_{ij} &= [r_1\cdots r_{j-1} i \mid c_1\cdots c_j] \in A
&&(j\le t,\ i\ge r_j)\\
v_{ij} &= \util_{ij}\drctilinv{j} \in \Arc &&(j\le t,\ i> r_j)\\
y_{ij} &= Y_{ij}Y_{r_jj}^{-1} \in \Rrplus \qquad \text{(as in (2.8))}
&&(j\le t,\ i> r_j). \endxalignat$$ 
(While the given expressions
for $v_{ij}$ and $y_{ij}$ also make sense for $j\le t$ and $i=r_j$,
they yield
$v_{r_jj} =1 \in \Arc$ and $y_{r_jj}=1 \in \Rrplus$. It is more
convenient for the proof to exclude these possibilities.) Recall that
the $y_{ij}$ generate
$\Brplus$.

{\bf Claim 1}: There exists a homomorphism $\phi^+: \Brplus
\rightarrow \Arc$ such that $\phi^+( y_{ij} )= v_{ij}$ for $j\le
t$ and $i>r_j$.

To prove this, we must show that the $v_{ij}$ satisfy the analogs of
the relations (2.8)(1), i.e., the corresponding equations with all
$y$'s replaced by $v$'s. We first check that the
$u_{ij}$ satisfy the relations (1) below. 
The first relation follows from Lemma 2.9(a); for the second, observe
that $\{r_1,\dots,r_{j-1},i\} \subset
\{r_1,\dots,r_{m-1},l\}$ and $\{c_1,\dots,c_j\} \subset
\{c_1,\dots,c_m\}$ in that case.
$$u_{lm}u_{ij} = \cases q^{-1}u_{ij}u_{lm} &\qquad (i<l,\, j\ge m)\\
  u_{ij}u_{lm} &\qquad (i\le l,\, j<m,\, i\in \{r_j,\dots,
r_{m-1},l\}). \endcases \tag1$$ 
Furthermore, when $i<l$ and $j<m$
but $i \notin \{r_j,\dots, r_{m-1}\}$, Lemma 2.9(b) implies that
$$u_{ij}u_{lm} -q^{-1} u_{lm}u_{ij} \congmodKrc \cases \qhat
u_{lj}u_{im} &\qquad (i\ge r_m)\\ 0 &\qquad (i<r_m). \endcases \tag2$$
Since $u_{r_ll}= \drc{l}$ for $l\le t$, the relations (1) and (2)
yield commutation relations for the $\drc{l}$ and $u_{ij}$, which we
combine in the following form:
$$\aligned \drc{l}u_{ij} &=\ u_{ij}\drc{l} \qquad\qquad \cases
(l<j), \text{\ or} &\\
 (l\ge j \text{\ and\ } i\in \{r_j,\dots,r_l\}) & \endcases \\
\drc{l}u_{ij}  &\congmodKrc  qu_{ij}\drc{l} \quad\qquad (l\ge j
\text{\ and\ } i\notin \{r_j,\dots,r_l\}). \endaligned \tag3$$ 
(Note that when $l<j$, we have $i\ge r_j> r_l$.)

It follows from (1), (2), and (3) that the elements $v_{ij} \in \Arc$
indeed satisfy the analogs of the relations (2.8)(1), as desired. This
establishes Claim 1.

Next, set 
$$\xalignat2 w_{lm} &= [r_1\cdots r_{l} \mid c_1\cdots c_{l-1} m]
\in A &&(l\le t,\ m\ge c_l)\\
t_{lm} &= \wtil_{lm}\drctilinv{l} \in \Arc &&(l\le t,\ m> c_l)\\
z_{lm} &= Z_{lm}Z_{lc_l}^{-1} \in \Rcminus &&(l\le t,\ m> c_l),
\endxalignat$$ 
and recall that the
$z_{lm}$ generate $\Bcminus$. 

{\bf Claim 2}. There exists a homomorphism $\phi^-: \Bcminus
\rightarrow \Arc$ such that $\phi^-(z_{lm})= t_{lm}$ for $l\le t$ and
$m>c_l$.

This follows from Claim 1 by symmetry.

{\bf Claim 3}. Each $v_{ij}$ commutes with each $t_{lm}$.

We first collect the following commutation relations between the
$u_{ij}$ and the
$w_{lm}$:
$$\aligned u_{ij}w_{lm} &= \ w_{lm}u_{ij} \quad\qquad\qquad \cases
(j=l), \text{\ or} &\\
 (j<l \text{\ and\ } i\in \{r_{j+1},\dots,r_l\}), \text{\ or} &\\
  (j>l \text{\ and\ } m\in \{c_{l+1},\dots,c_j\}) & \endcases \\
u_{ij}w_{lm} &\congmodKrc \cases q^{-1}w_{lm}u_{ij} &\quad (j<l
\text{\ and\ } i\notin \{r_{j+1},\dots,r_l\})\\
  qw_{lm}u_{ij} &\quad (j>l \text{\ and\ } m\notin
\{c_{l+1},\dots,c_j\}). \endcases \endaligned \tag4$$ 
The first equation in (4) follows from Lemma 5.3(c). The next two
equations hold because $\{r_1,\dots,r_{j-1},i\}
\subset \{r_1,\dots,r_l\}$ and $\{c_1,\dots,c_j\} \subset
\{c_1,\dots,c_{l-1},m\}$ in the first case, while those inclusions
are reversed in the second case. Finally, the last two relations
follow from Lemma 2.9 and Corollary 2.10. 

Commutation relations between the $\drc{l}$ and the $u_{ij}$ are given
in (3), and, by symmetry, we also have
$$\aligned \drc{j}w_{lm} &= \ w_{lm}\drc{j} \qquad\qquad \cases (j<l),
\text{\ or} &\\
 (j\ge l \text{\ and\ } m\in
\{c_l,\dots,c_j\})& \endcases \\
\drc{j}w_{lm} &\congmodKrc  qw_{lm}\drc{j} \quad\qquad (j\ge l
\text{\ and\ } m\notin \{c_l,\dots,c_j\}). \endaligned \tag5$$  
It follows from
(3), (4), and (5) that $v_{ij}$ indeed commutes with $t_{lm}$,
and Claim 3 is proved.

Combining Claims 1, 2, and 3, we see that there exists a homomorphism
$$\phi: \Brplus \otimes \Bcminus \rightarrow \Arc$$
such that $\phi( y_{ij}\otimes 1 )= v_{ij}$ and $\phi(1\otimes z_{lm})=
t_{lm}$ for all $i,j,l,m$. 

{\bf Claim 4}. $\phi$ extends to a homomorphism
$\Brc \rightarrow \Arc$ such that $\phi( Y_{r_ss} \otimes Z_{sc_s} )= 
\drctil{s}\drctilinv{s-1}$ for $s\le t$.

In view of the relations given in (2.1), we see that
$$\aligned Y_{r_ss}y_{ij}Y_{r_ss}^{-1} &= \cases q^{-1}y_{ij} &\quad
(i=r_s)\\ qy_{ij} &\quad (j=s)\\ y_{ij} &\quad (i\ne r_s,\ j\ne s)
\endcases\\
Z_{sc_s}z_{lm}Z_{sc_s}^{-1} &= \cases qz_{lm} &\quad (l=s)\\
q^{-1}z_{lm} &\quad (m=c_s)\\ z_{lm} &\quad (l\ne s,\ m\ne c_s).
\endcases \endaligned \tag6$$
On the other hand, it follows from (3) and (5) that the units
$\drctil{s}\drctilinv{s-1}$ in $\Arc$ normalize the elements $v_{ij}$
and $t_{lm}$ in exactly the same way, that is,
$$\aligned \drctil{s}\drctilinv{s-1}v_{ij} \drctil{s-1}\drctilinv{s} &=
\cases q^{-1}v_{ij} &\quad (i=r_s)\\ qv_{ij} &\quad (j=s)\\ v_{ij}
&\quad (i\ne r_s,\ j\ne s)
\endcases\\
\drctil{s}\drctilinv{s-1}t_{lm} \drctil{s-1}\drctilinv{s} &= \cases
qt_{lm} &\quad (l=s)\\ q^{-1}t_{lm} &\quad (m=c_s)\\ t_{lm} &\quad
(l\ne s,\ m\ne c_s).
\endcases \endaligned \tag7$$
Claim 4 follows from (6), (7), and Lemma 2.7.

{\bf Claim 5}. $\phi= \betatilrc^{-1}$. 

Since $\betatilrc$
is surjective (Lemma 2.5), it is enough to show that $\phi\betatilrc$
is the identity on $\Arc$. First note, using (2.4), that
$$\phi\betatilrc(\drctil{s}\drctilinv{s-1}) =\phi(Y_{r_ss} \otimes
Z_{sc_s})= \drctil{s}\drctilinv{s-1}$$
for all $s$, whence $\phi\betatilrc(\drctil{s})= \drctil{s}$ for all
$s$. Next, for $j\le t$ and $i>r_j$ we have
$$\align \phi\betatilrc( \util_{ij}) &= \phi\bigl( (Y_{r_11}\otimes
Z_{1c_1})\cdots (Y_{r_{j-1},j-1}\otimes Z_{j-1,c_{j-1}})
(Y_{ij}\otimes Z_{jc_j}) \bigr)\\
 &= \phi\bigl( (y_{ij} \otimes 1) (Y_{r_11}\otimes
Z_{1c_1})\cdots (Y_{r_jj}\otimes Z_{jc_j}) \bigr)\\
 &= \phi\bigl( (y_{ij} \otimes 1) \betarc(\drctil{j}) \bigr) = v_{ij}
\drctil{j} = \util_{ij}. \endalign$$ 
By symmetry,
$\phi\betatilrc(\wtil_{lm})= \wtil_{lm}$ for $l\le t$ and
$m>c_l$. Therefore $\phi\betatilrc$ at least equals the identity on
the subalgebra $C$ of $\Arc$ generated by (the image of) the set
$$\{ (\drctil{s})^{\pm1} \mid s\le t\} \cup \{ \util_{ij} \mid j\le t,\
i>r_j \} \cup \{ \wtil_{lm} \mid l\le t,\ m>c_l \}.$$

To finish the proof, we just need to show that $C= \Arc$, that is,
that $\Xtil_{ij} \in C$ for all $i,j$. This is clear in case $i<r_1$
or $j<c_1$, since in those cases $\Xtil_{ij} =0$. We also have
$\Xtil_{ic_1}=
\util_{i1} \in C$ for $i\ge r_1$ and $\Xtil_{r_1j}= \wtil_{1j} \in C$
for $j\ge c_1$. Hence, $\Xtil_{ij} \in C$ whenever $i\le r_1$ or $j\le
c_1$.

Now let $1<l\le t$ and assume that $\Xtil_{ij} \in C$ whenever $i\le
r_{l-1}$ or $j\le c_{l-1}$. For $r_{l-1}<i \le t$ and $c_{l-1}< j\le
c_l$, it follows from Corollary
5.5(b1) that
$$\multline (-q)^{l-1} [r_1\cdots r_{l-1}i \mid c_1\cdots c_{l-1}j]=
(-q)^{l-1} X_{ij} [r_1\cdots r_{l-1} \mid c_1\cdots c_{l-1}]\\ 
+ \sum_{s=1}^{l-1} (-q)^{s-1} X_{ic_s} [ r_1\cdots r_{l-1}\mid
c_1\cdots \widehat{c_s}\cdots c_{l-1}j ]. \endmultline \tag8$$ 
Most of
the cosets of the factors in (8) can be seen to lie in
$C$ right away. For instance, the coset of the left hand side is zero
if either
$i<r_l$ or $j<c_l$, and it equals $\util_{il}$ if $i\ge r_l$ and
$j=c_l$. For
$s<l$, we have $\Xtil_{ic_s} \in C$ by the inductive hypothesis, and
similarly the coset $[r_1\cdots r_{l-1}\mid
c_1\cdots \widehat{c_s}\cdots c_{l-1}j] +\Krc$ is in $C$ since it is a
linear combination of products of elements $\Xtil_{ab}$ with $a\le
r_{l-1}$. Consequently, we obtain from (8) that 
$$X_{ij} [r_1\cdots r_{l-1} \mid c_1\cdots c_{l-1}] +\Krc \in C.$$
In other words, $\Xtil_{ij} \drctil{l-1} \in C$, whence
$\Xtil_{ij} \in C$. Thus, $\Xtil_{ij} \in C$ whenever $j\le c_l$. By
symmetry, $\Xtil_{ij} \in C$ whenever $i\le r_l$. 

The above induction proves that $\Xtil_{ij} \in C$ whenever $i\le r_t$
or $j\le c_t$. If there exist indices $i>r_t$ and $j>c_t$, we have
$$\multline (-q)^{t} [r_1\cdots r_ti \mid c_1\cdots c_tj]=
(-q)^t X_{ij} [r_1\cdots r_t \mid c_1\cdots c_t]\\
+\sum_{s=1}^t (-q)^{s-1} X_{ic_s} [ r_1\cdots r_t\mid
c_1\cdots \widehat{c_s}\cdots c_tj ] \endmultline \tag9$$
by Corollary 5.5(b1), from which we see as above that $\Xtil_{ij} \in C$.
(Note that the left hand side of (9) necessarily lies in $\Krc$
because it involves a $(t+1)\times (t+1)$ quantum minor.) Therefore all
$\Xtil_{ij}
\in C$, and the proof is complete.
\qed\enddemo

\head 3. Tensor product decompositions of $H$-primes \endhead

Throughout this section, we assume that $q$ is not a root of unity; we
shall place reminders of this hypothesis in the relevant results. Thus,
by \cite{\GLet, Theorem 3.2}, all primes of
$A$ are completely prime. Since this property survives in factors and
localizations, all primes in the algebras
$\Arc$, $\Rrplus$, and $\Rcminus$ are completely prime, and also in
$\Brc$ because of Theorem 2.11. We have already observed that the
algebras $\Rrpluszero$ and $\Rcminuszero$ are iterated skew polynomial
algebras over $k$, and so is their tensor product. The iterated skew
polynomial structure of $\Rrpluszero \otimes \Rcminuszero$ is easily
seen to satisfy the hypotheses of \cite{\GLet, Theorem 2.3}, and
thus all its primes are completely prime. Consequently, all
primes in the localizations $\ror$ are completely
prime.

\definition{3.1} In order to deal with $H$-primes in tensor products,
we need the following rationality property. Suppose that $S$ is a
noetherian $k$-algebra and that $G$ is a group acting on $S$ by
$k$-algebra automorphisms. We say that a $G$-prime $P$ of $S$ is
{\it strongly $G$-rational\/} provided the algebra $Z(\fract S/P)^G$,
the fixed ring of the center of the Goldie quotient ring of
$S/P$ under the induced $G$-action, equals $k$.

 By \cite{\specstrat, (5.7)(i)} (cf\. \cite{\barcelona, Theorem
II.5.14}) and \cite{\barcelona, Corollary II.6.5},
$A$ has only finitely many $H$-primes, and they are all completely
prime and strongly $H$-rational. These properties carry over to $\Arc$,
$\Rrplus$, and $\Rcminus$. Analogous results \cite{\barcelona, Theorems
II.5.12 and II.6.4} imply that
$\Rrpluszero \otimes \Rcminuszero$ has only finitely many $(H\times
H)$-primes, and they are all completely prime and strongly $(H\times
H)$-rational. These properties now carry over to $\ror$.

Identify $H$ with the subgroup
$$\Htil= \bigl( (\kx)^n \times \{1\}^n \bigr) \times \bigl( \{1\}^n
\times (\kx)^n \bigr) \vartriangleleft H\times H$$
in the obvious way. With this identification,
$\betarc$ and $\betatilrc$ are $H$-equivariant. In particular, it
follows that
$\Brc$ has only finitely many
$\Htil$-primes, and they are all completely prime and strongly
$\Htil$-rational.
\enddefinition

\proclaim{3.2\. Lemma} For $i=1,2$, let $A_i$ be a $k$-algebra, $H_i$
a group acting on $A_i$ by $k$-algebra automorphisms, and $P_i$ an
$H_i$-prime ideal of $A_i$. Set $P= (P_1\otimes A_2) +(A_1\otimes
P_2)$, and let $H_1\times H_2$ act on $A_1\otimes A_2$ in the natural
manner.

{\rm (a)} If each $A_i/P_i$ is $H_i$-simple and $Z(A_1/P_1)^{H_1}=k$,
then
$(A_1\otimes A_2)/P$ is $(H_1\times H_2)$-simple.

{\rm(b)} If each $A_i$ is noetherian and $P_1$ is strongly
$H_1$-rational, then
$P$ is an $(H_1\times H_2)$-prime ideal of $A_1\otimes A_2$. Moreover,
$P$ is the only $(H_1\times H_2)$-prime ideal of $A_1\otimes A_2$ that
contracts to
$P_1\otimes 1$ in $A_1\otimes 1$ and to $1\otimes P_2$ in
$1\otimes A_2$.

{\rm (c)} If each $A_i$ is noetherian and each $P_i$ is strongly
$H_i$-rational, then $P$ is strongly $(H_1\times H_2)$-rational.
\endproclaim

\demo{Proof} Since $(A_1\otimes A_2)/P \cong (A_1/P_1)\otimes
(A_2/P_2)$, there is no loss of generality in assuming that each
$P_i=0$.

(a) This is a standard shortest length argument. Let $I$ be a nonzero
$(H_1\times H_2)$-ideal of $A_1\otimes A_2$, and let $m$ be the
shortest length for nonzero elements of $I$ (as sums of pure tensors).
Choose a nonzero element
$$x= b_1\otimes
c_1 +\dots+ b_m\otimes c_m \in I$$
of length $m$, where the $b_j\in A_1$ and $c_j\in A_2$,
and note that the $c_j$ are linearly independent over $k$. Now the set
$$\{ b\in A_1\mid (b\otimes c_1+ A_1\otimes c_2 +\dots+ A_1\otimes
c_m) \cap I\ne \varnothing \}$$ 
is a nonzero $H_1$-ideal of $A_1$, and so
it equals $A_1$. Thus, without loss of generality, $b_1=1$. For any
$a\in A_1$, we now have
$x(a\otimes 1)- (a\otimes 1)x$ in $I$ with length less than
$m$, whence $x(a\otimes 1)- (a\otimes 1)x=0$ and so $b_ja=ab_j$ for
all $j$. For any $h\in H_1$, we have $(h,1)(x)-x$ in $I$ with length
less than $m$, whence
$(h,1)(x)-x=0$ and so $h(b_j)=b_j$ for all $j$. Therefore all
$b_j$ lie in $Z(A_1)^{H_1}=k$. It follows that $x\in 1\otimes A_2$,
whence $m=1$ and
$x= 1\otimes c_1$. Consequently, the set $\{ c\in A_2\mid 1\otimes
c\in I\}$ is a nonzero $H_2$-ideal of $A_2$, and so it equals $A_2$.
Therefore $1\otimes 1\in I$, proving that $I= A_1\otimes A_2$.

(b) Each $A_i$ is an $H_i$-prime noetherian ring, and so is semiprime.
Let
$\C_i$ be the set of regular elements in $A_i$, and note that the set
$$\C= \{ c_1\otimes c_2\mid c_i\in \C_i \}$$ 
is an $(H_1\times
H_2)$-stable denominator set in
$A_1\otimes A_2$, consisting of regular elements. Each $A_i\C_i^{-1}$
is
$H_i$-simple artinian, and $Z(A_1\C_1^{-1})^{H_1} =k$ by our
hypothesis on $P_1$. Thus by part (a), the localization
$$(A_1\otimes A_2)\C^{-1}= A_1\C_1^{-1} \otimes A_2\C_2^{-1}$$ is
$(H_1\times H_2)$-simple. It follows that each nonzero $(H_1\times
H_2)$-ideal of $A_1\otimes A_2$ meets $\C$; in particular, $A_1\otimes
A_2$ is an $(H_1\times H_2)$-prime ring.

Now let $Q$ be any $(H_1\times H_2)$-prime ideal of $A_1\otimes A_2$
that contracts to zero in both $A_1\otimes 1$ and $1\otimes A_2$. Then
$Q$ is a semiprime ideal, disjoint from both $\C_1\otimes 1$ and
$1\otimes \C_2$. If some prime ideal $Q_0$ minimal over $Q$ meets
$\C_1\otimes 1$, then $h(Q_0)$ meets $\C_1\otimes 1$ for all $h\in
H_1\times H_2$. But since $Q$ is a finite intersection of some of the
$h(Q_0)$, it would follow that $Q$ meets
$\C_1\otimes 1$, a contradiction. Therefore $\C_1\otimes 1$ is
disjoint from all primes minimal over $Q$, whence $\C_1\otimes 1$ is
regular modulo $Q$. Likewise, $1\otimes \C_2$ is regular modulo $Q$.
It follows that $\C$ is disjoint from $Q$, and therefore $Q=0$.

(c) After localization, we can assume that each $A_i$ is $H_i$-simple
artinian. By part (a), $A_1\otimes A_2$ is now $(H_1\times
H_2)$-simple. Consider an element $u$ in
$Z(\fract (A_1\otimes A_2))^{H_1\times H_2}$. The set
$\{a\in A_1\otimes A_2\mid au\in A_1\otimes A_2\}$ is a nonzero
$(H_1\times H_2)$-ideal of $A_1\otimes A_2$, and so it equals
$A_1\otimes A_2$. Therefore
$u\in A_1\otimes A_2$. Now write $u= v_1\otimes w_1 +\dots+ v_t\otimes
w_t$ for some $v_j\in A_1$ and some linearly independent $w_j\in A_2$.
Since $u$ is fixed by $H_1\times 1$ and commutes with $A_1\otimes 1$,
we see that all
$v_j\in Z(A_1)^{H_1}=k$. Hence, $u=1\otimes w$ for some $w\in A_2$.
But then
$w\in Z(A_2)^{H_2}=k$, and therefore $u\in k$. \qed\enddemo

\proclaim{3.3\. Proposition} For $i=1,2$, let $A_i$ be a noetherian
$k$-algebra and
$H_i$ a group acting on $A_i$ by $k$-algebra automorphisms. Assume
that all
$H_1$-primes of $A_1$ are strongly $H_1$-rational. Then the rule
$(P_1,P_2)
\mapsto (P_1\otimes A_2) +(A_1\otimes P_2)$ provides a bijection
$$(H_1\operatorname{-spec} A_1)\times (H_2\operatorname{-spec} A_2)
\longrightarrow (H_1\times H_2)\operatorname{-spec} (A_1\otimes A_2).$$
\endproclaim

\demo{Proof} Lemma 3.2(b) shows that the given rule maps
$(H_1\operatorname{-spec} A_1)\times (H_2\operatorname{-spec} A_2)$ to
$(H_1\times H_2)\operatorname{-spec} (A_1\otimes A_2)$.

 Now consider an $(H_1\times H_2)$-prime $P$ in $A_1\otimes A_2$. Let
$P_1$ and
$P_2$ be the inverse images of $P$ under the natural maps $A_i
\rightarrow A_1\otimes A_2$. Then each $P_i$ is an $H_i$-ideal of
$A_i$, and $(P_1\otimes A_2) +(A_1\otimes P_2) \subseteq P$. There are
$H_1$-primes $Q_1,\dots,Q_t$ in $A_1$, containing $P_1$, such that
$Q_1Q_2
\cdots Q_t\subseteq P_1$. Then the
$Q_i\otimes A_2$ are $(H_1\times H_2)$-ideals of $A_1\otimes A_2$ such
that
$$(Q_1\otimes A_2)(Q_2\otimes A_2) \cdots(Q_t\otimes A_2) \subseteq
P_1\otimes A_2 \subseteq P.$$ Consequently, some $Q_j\otimes
A_2\subseteq P$, whence $Q_j\subseteq P_1$, and so $Q_j=P_1$. This
shows that $P_1$ is $H_1$-prime. Similarly, $P_2$ is
$H_2$-prime.

By Lemma 3.2(b), $(P_1\otimes A_2) +(A_1\otimes P_2)$ is an $(H_1\times
H_2)$-prime of $A_1\otimes A_2$, and it is the only $(H_1\times
H_2)$-prime of $A_1\otimes A_2$ that contracts to $P_1\otimes 1$ in
$A_1\otimes 1$ and to $1\otimes P_2$ in $1\otimes A_2$. Therefore $P=
(P_1\otimes A_2) +(A_1\otimes P_2)$. It is clear that $P_1$ and $P_2$
are unique, since $P_i$ equals the inverse image of $(P_1\otimes A_2)
+(A_1\otimes P_2)$ under the natural map $A_i \rightarrow A_1\otimes
A_2$.
\qed\enddemo

\proclaim{3.4\. Lemma} {\rm [$q$ not a root of unity]} Let $(\rc) \in
\RC$. If $\Ptil$ is an $\Htil$-prime of $\Brc$, then there exists a
unique $(H\times H)$-prime $Q$ of $\ror$ such that $Q\cap \Brc=
\Ptil$. \endproclaim

\demo{Proof} Since $\Brc$ has only finitely many $\Htil$-primes and
$H\times H$ just permutes them, the $(H\times H)$-orbit of $\Ptil$ in
$\spec \Brc$ is finite. Since $\Ptil$ is prime, it now follows from
\cite{\barcelona, Proposition II.2.9} that $\Ptil$ must be invariant
under $H\times H$. In view of Lemma 2.6, $Q= \Ptil (\ror)$ is an
$(H\times H)$-in\-var\-i\-ant prime of $\ror$ such that $Q\cap \Brc= \Ptil$.
It remains to show that if $Q'$ is any $(H\times H)$-prime of
$\ror$ that contracts to $\Ptil$, then $Q'=Q$. Note that $Q' \supseteq
Q$ by definition of $Q$.

Set $G= (\kx)^n$, let $\phi: G\rightarrow \bigl( \{1\}^n\times (\kx)^n
\bigr) \times \bigl( (\kx)^n\times \{1\}^n \bigr) \subset
H\times H$ be the homomorphism given by the rule
$$\phi(\alpha_1,\dots,\alpha_n)=
(1,\dots,1,\alpha_1^{-1},\dots,\alpha_n^{-1}, \alpha_1,\dots,\alpha_n,
1,\dots,1),$$
and use $\phi$ to pull back the action of $H\times H$ on $\ror$ to an
action of $G$. With respect to this $G$-action, $\Brc$ is
generated by fixed elements, and each of the elements $1\otimes
Z_{sc_s}$ is a $G$-eigenvector with eigenvalue equal to the projection
$(\alpha_1,\dots,\alpha_n) \mapsto \alpha_s$. In view of Lemma 2.6 and
the fact that
$k$ is infinite, it follows that the $G$-eigenspaces of $\ror$ are the
subspaces $\Brc \bigl( 1\otimes (Z_{1c_1}^{m_1} Z_{2c_2}^{m_2} \cdots
Z_{tc_1}^{m_t}) \bigr)$ for $(m_1,\dots,m_t)\in \ZZ^t$. Consequently,
any $G$-eigenvector in $\ror$ has the form $du$ where $d\in
\Brc$ and $u$ is a unit. If $du\in Q'$, then $d\in Q'\cap \Brc= \Ptil$,
whence $du\in Q$. Since $Q'$ is $G$-in\-var\-i\-ant, we conclude that
$Q'=Q$, as desired. \qed\enddemo

\definition{3.5} Set $\Hspec_{\rc} A= (\Hspec A)\cap (\spec_{\rc} A)$
for $(\rc)\in \RC$. These sets partition $\Hspec A$ because of
Corollary 1.10.
\enddefinition

\proclaim{Theorem} {\rm [$q$ not a root of unity]} For each  $(\rc)
\in \RC$, there is a bijection
$$(\Hspec \Rrplus) \times (\Hspec \Rcminus) \longrightarrow 
\Hspec_{\rc} A$$
given by the rule $(Q^+,Q^-) \mapsto \betarc^{-1} \bigl(
(Q^+\otimes
\Rcminus) +(\Rrplus\otimes Q^-) \bigr)$.
\endproclaim

\demo{Proof} If $Q^+\in \Hspec \Rrplus$ and $Q^- \in \Hspec
\Rcminus$, then Proposition 3.3 shows that the ideal $Q= (Q^+\otimes
\Rcminus) +(\Rrplus\otimes Q^-)$ is an $(H\times H)$-prime of $\ror$.
In particular, $Q$ is completely prime, and so $Q\cap \Brc$ is an
$\Htil$-prime of $\Brc$, whence $\betarc^{-1}(Q)= \betarc^{-1}(Q\cap
\Brc)$ is an $H$-prime of $A$ lying in $\spec_{\rc} A$. This shows
that the given rule does define a map from $(\Hspec \Rrplus) \times
(\Hspec \Rcminus)$ to $\Hspec_{\rc} A$.

Now consider an arbitrary $H$-prime $P$ in $\Hspec_{\rc} A$. Then $P$
induces an
$H$-prime in $\Arc$ that contracts to $P$ under the localization map.
In view of Theorem 2.11, it follows that
$\betarc(P)$ induces an
$\Htil$-prime $\Ptil$ of $\Brc$ such that
$\betarc^{-1}(\Ptil) =P$. By Lemma 3.4, there is a unique $(H\times
H)$-prime $Q$ of $\ror$ such that $Q\cap \Brc= \Ptil$. Then
Proposition 3.3 implies that $Q= (Q^+\otimes
\Rcminus) +(\Rrplus\otimes Q^-)$ for some $H$-primes $Q^+$ in
$\Rrplus$ and $Q^-$ in $\Rcminus$. Thus,
$$\betarc^{-1} \bigl( (Q^+\otimes \Rcminus) +(\Rrplus\otimes Q^-)
\bigr)= \betarc^{-1}(Q)= \betarc^{-1}(Q\cap \Brc)= \betarc^{-1}(\Ptil)
=P.$$

It remains to show that $Q^+$ and $Q^-$ are unique.
Consider any $T^+\in \Hspec \Rrplus$ and $T^-\in \Hspec
\Rcminus$ such that $P= \betarc^{-1}(T)$ where
$T= (T^+\otimes \Rcminus) +(\Rrplus\otimes T^-)$. As in the first
paragraph of the proof,
$T$ is an
$(H\times H)$-prime of $\ror$ and
$T\cap \Brc$ is an $\Htil$-prime of $\Brc$. Since $\betarc^{-1}(T\cap
\Brc)= P$, we must have $T\cap \Brc=
\Ptil$, whence $T=Q$ by the uniqueness of $Q$. Therefore Proposition
3.3 shows that
$T^+=Q^+$ and
$T^-=Q^-$, as desired. \qed\enddemo

\definition{3.6} Fix $t\in \{0,1,\dots,n\}$, and let $\Hspec^{[t]} A$
be the set of those $H$-primes of $A$ which contain all $(t+1)\times
(t+1)$ quantum minors but not all $t\times t$ quantum minors. By
Corollary 1.10,
$$\Hspec^{[t]} A= \bigsqcup_{(\rc)\in \RC_t} \Hspec_{\rc}A,$$
and consequently Theorem 3.5 implies that
$$|\Hspec^{[t]} A|= \sum_{(\rc)\in\RC_t} |\Hspec\Rrplus| \cdot
|\Hspec\Rcminus|.$$
If $\Rt$ denotes the set of sequences $(r_1,\dots,r_t) \in \NN^t$ with
$1\le r_1< \cdots< r_t\le n$, then $\RC_t= \Rt\times\Rt$. For each
$\bfr\in\Rt$, the automorphism $\bftau$ of $A$ discussed in (5.1)
induces an isomorphism $\overline{\bftau} : \Rrplus
\rightarrow R^-_{\bfr}$. While $\overline{\bftau}$ is not
$H$-equivariant, there is an automorphism $\gamma$ of $H$, given by
$(\alpha_1,\dots,\alpha_n, \beta_1,\dots,\beta_n) \mapsto
(\beta_1,\dots,\beta_n, \alpha_1,\dots,\alpha_n)$, such that the
following diagram commutes:
$$\CD H @>{\gamma}>> H\\
@VVV @VVV\\
\Aut \Rrplus @>{\overline{\bftau}^*}>> \Aut R^-_{\bfr} \endCD$$
(Here the vertical arrows denote the standard actions of $H$ on
$\Rrplus$ and $R^-_{\bfr}$.) Hence, $\overline{\bftau}$ provides a
bijection of $\Hspec \Rrplus$ onto $\Hspec R^-_{\bfr}$. Therefore
$$|\Hspec^{[t]} A|= \sum_{\bfr,\bfc\in\Rt} |\Hspec\Rrplus| \cdot
|\Hspec R^+_{\bfc}|= \biggl(\; \sum_{\bfr\in\Rt} |\Hspec\Rrplus|
\biggr)^2,$$ a perfect square. These numbers are known in three cases:
$$|\Hspec^{[t]} A|= \cases 1 &\qquad (t=0)\\ (2^n-1)^2 &\qquad (t=1)\\
(n!)^2 &\qquad (t=n). \endcases$$
The case when $t=0$ is trivial, and the case when $t=1$ is given by
\cite{\GLmurcia, Corollary 3.5}. For the remaining case, note first
that $\Hspec^{[n]} A \approx \Hspec \Oq(GL_n(k))$. It can be checked
that there is a bijection between $\Hspec \Oq(GL_n(k))$ and the set of
winding-in\-var\-i\-ant primes of $\Oq(SL_n(k))$ (e.g., see \cite{\barcelona,
Lemma II.5.16}), and it follows from the work of Hodges and Levasseur
\cite{\HoLetwo} that the latter set is in bijection with the double
Weyl group $S_n\times S_n$ (cf\. \cite{\barcelona, Corollary II.4.12}).
\enddefinition

\definition{3.7} As a corollary of Theorem 3.5, we obtain the
following less specific but more digestible result.
\enddefinition

\proclaim{Corollary} {\rm [$q$ not a root of unity]} Set
$$R^+= A/\langle X_{ij}\mid i<j \rangle \qquad\qquad \text{and}
\qquad\qquad R^-= A/\langle X_{ij} \mid i>j \rangle,$$
 let $\pi^{\pm} : A\rightarrow R^{\pm}$ denote the quotient maps, and
let $\beta$ denote the composition
$$A@>{\ \Delta\ }>> A\otimes A @>{\pi^+ \otimes \pi^-}>> R^+
\otimes R^-.$$
Given any $H$-prime $P$ in $A$, there exist $H$-primes $P^{\pm}$ in
$R^{\pm}$ such that
$$P= \beta^{-1} \bigl( (P^+\otimes R^-) +(R^+\otimes P^-) \bigr).$$
\endproclaim

\demo{Proof} By Corollary 1.10 and Theorem 3.5,
$P= \betarc^{-1} \bigl( (Q^+\otimes \Rcminus)
+(\Rrplus\otimes Q^-) \bigr)$
for some $(\rc)\in \RC$ and some $H$-primes $Q^+$ in $\Rrplus$ and
$Q^-$ in $\Rcminus$. Observe that $X_{ij} \in \ker \pirpluszero$ when
$i<j$, and that $X_{ij} \in \ker \picminuszero$ when $i>j$. Hence,
there are surjective $k$-algebra homomorphisms $\tau^+ : R^+
\rightarrow
\Rrpluszero$ and $\tau^- : R^- \rightarrow \Rcminuszero$ such that
$\tau^+\pi^+= \pirpluszero$ and $\tau^-\pi^-= \picminuszero$.
Consequently, $(\tau^+ \otimes \tau^-)\beta =\betarc$ (with the
obvious adjustment of codomains).

Next, observe that $Q^+_0= Q^+\cap \Rrpluszero$ and $Q^-_0= Q^- \cap
\Rcminuszero$ are $H$-primes of $\Rrpluszero$ and $\Rcminuszero$,
respectively, whence the ideals $P^{\pm}= (\tau^{\pm})^{-1}(Q^{\pm}_0)$
are $H$-primes in $R^{\pm}$. Finally,
$$P= \bigl( (\tau^+ \otimes \tau^-)\beta \bigr)^{-1} \bigl( (Q^+_0
\otimes \Rcminuszero) +(\Rrpluszero \otimes Q^-_0) \bigr)= \beta^{-1}
\bigl( (P^+\otimes R^-) +(R^+\otimes P^-) \bigr),$$
as desired. \qed\enddemo

\head 4. Illustration: $\OqMtwo$ \endhead

Theorem 3.5 opens a potential route to computing the $H$-primes of
$A$ in the generic case: If we can find all
the $H$-primes in each
$\Rrplus$ and
$\Rcminus$, we immediately obtain descriptions of all the $H$-primes
in $A$. Since these descriptions would be in terms of pullbacks of
$H$-primes from the algebras $\Rrplus \otimes \Rcminus$, it would
still remain to find generating sets for these ideals.

To illustrate the procedure, we sketch
the case where $n=2$, for which $\Hspec A$ is already known. In
\cite{\GLenseq}, we use the above process to compute
$\Hspec A$ when $n=3$.

\definition{4.1} Assume that $q$ is not a root of unity, and fix
$n=2$. There are only four choices for $\bfr$ and $\bfc$, namely
$\varnothing$, $(1)$, $(2)$, and $(1,2)$. The corresponding algebras
$\Rrplus$ and $\Rcminus$ are
$$\xalignat2 R^+_{\varnothing} &= A/\langle X_{11}, X_{12}, X_{21},
X_{22} \rangle =k  &R^-_{\varnothing} &=k\\
R^+_{(1)} &= \bigl( A/\langle X_{12}, X_{22} \rangle
\bigr)[X_{11}^{-1}]= k\langle Y_{11}^{\pm1}, Y_{21} \rangle 
&R^-_{(1)} &= k\langle Z_{11}^{\pm1}, Z_{12} \rangle\\
R^+_{(2)} &= \bigl( A/\langle X_{11}, X_{12}, X_{22} \rangle
\bigr)[X_{21}^{-1}]= k[Y_{21}^{\pm1}]  &R^-_{(2)} &= k[Z_{12}^{\pm1}]\\
R^+_{(1,2)} &= \bigl( A/ \langle X_{12}\rangle \bigr) [X_{11}^{-1},
X_{22}^{-1}]= k\langle Y_{11}^{\pm1}, Y_{21}, Y_{22}^{\pm1} \rangle 
&R^-_{(1,2)} &= k\langle Z_{11}^{\pm1}, Z_{12}, Z_{22}^{\pm1} \rangle.
\endxalignat$$
The $H$-primes in these algebras are easily computed:
$$\xalignat2 \Hspec R^+_{\varnothing} &= \bigl\{ \langle 0\rangle
\bigr\}  &\Hspec R^-_{\varnothing} &= \bigl\{ \langle 0\rangle
\bigr\}\\
\Hspec R^+_{(1)} &= \bigl\{ \langle 0\rangle,\ \langle Y_{21}\rangle
\bigr\}  &\Hspec R^-_{(1)} &= \bigl\{ \langle 0\rangle,\ \langle
Z_{12}\rangle \bigr\}\\
\Hspec R^+_{(2)} &= \bigl\{ \langle 0\rangle \bigr\}  &\Hspec
R^-_{(2)} &= \bigl\{ \langle 0\rangle \bigr\}\\
\Hspec R^+_{(1,2)} &= \bigl\{ \langle 0\rangle,\ \langle Y_{21}\rangle
\bigr\}  &\Hspec R^-_{(1,2)} &= \bigl\{ \langle 0\rangle,\ \langle
Z_{12}\rangle \bigr\}. \endxalignat$$

The only choice for $(\rc) \in \RC_0$ is $\bfr= \bfc= \varnothing$. In
this case, the only $H$-primes in $\Rrplus$ and $\Rcminus$ are the zero
ideals, and $\beta_{\varnothing,\varnothing}^{-1}( \langle 0\rangle)=
\langle X_{11}, X_{12}, X_{21}, X_{22} \rangle$, the augmentation
ideal of
$A$. We record this
$H$-prime using the symbol

\ignore
$$\xymatrixrowsep{0.1pc}\xymatrixcolsep{0.1pc}
\xymatrix{
\plb &\plb\\ \plb &\plb
}$$
\endignore
\medskip

\noindent to denote the generating set $\{ X_{11}, X_{12},
X_{21}, X_{22} \}$, the bullet in position $(i,j)$ being a marker
for the element $X_{ij}$.

Corresponding to the four pairs $(\rc) \in \RC_1$, there are nine
$H$-primes in $A$ of the form
$$\betarc^{-1} \bigl( (Q^+\otimes \Rcminus)
+(\Rrplus\otimes Q^-) \bigr)$$
where $Q^+$ is an $H$-prime in $\Rrplus$ and $Q^-$ is an $H$-prime in
$\Rcminus$. We can record generating sets for these ideals as follows,
continuing the notation introduced in the previous paragraph; here
$\circ$ is a placeholder and $\square$ denotes the $2\times2$ quantum
determinant.

\ignore
$$\xymatrixrowsep{0.1pc}\xymatrixcolsep{0.1pc}
\xymatrix{
 &&&&\place \edge[15,0] &&&& &&&&\place \edge[15,0] \\
 \\   \\
\place \edge[0,16] &&&&&&\dropup{3}{\langle 0\rangle}
\dropup{8}{R^-_{(1)}} &&&\dropup{3}{\langle Z_{12}\rangle}
&&&&&\dropup{3}{\langle 0\rangle} \dropup{8}{R^-_{(2)}}
&&&\place \\  \\
 &&&\dropleft{3}{\langle 0\rangle} \dropleft{9}{R^+_{(1)}}
&&&\hrzvrt &&\plc &\plb &&&&\plb &\plc \\
 &&&&&&\hrz &&\plc &\plb &&&&\plb &\plc \\   \\
 &&&\dropleft{3}{\langle Y_{21}\rangle} &&&\plc &\plc
&&\plc &\plb &&&&\plb &\plc \\
 &&&&&&\plb &\plb &&\plb &\plb &&&&\plb &\plb \\   \\
\place \edge[0,16] &&&&& &&&&& &&&&&&\place \\   \\
 &&&\dropleft{3}{\langle 0\rangle} \dropleft{9}{R^+_{(2)}}
&&&\plb &\plb &&\plb &\plb &&&&\plb &\plb \\
 &&&&\place &&\plc &\plc &&\plc &\plb &&\place &&\plb
&\plc \\
&&&&\place &&&& &&&&\place
}$$
\endignore
\medskip

\noindent (See, e.g., \cite{\GLmurcia, Theorem 1.1} for a proof that
the quantum determinant generates the kernel of $\beta_{(1),(1)}$. We
leave it to the reader to check that the other $H$-primes are
generated as indicated.)

Finally, the only choice for $(\rc) \in \RC_2$ is $\bfr= \bfc= (1,2)$,
and there are four 
$H$-primes in $A$ of the form
$$\beta_{(1,2),(1,2)}^{-1} \bigl( (Q^+\otimes R^-_{(1,2)})
+(R^+_{(1,2)} \otimes Q^-) \bigr)$$
with $Q^{\pm} \in \Hspec R^{\pm}_{(1,2)}$. We record generating sets
for these ideals as follows:

\ignore
$$\xymatrixrowsep{0.1pc}\xymatrixcolsep{0.1pc}
\xymatrix{
 &&&&\place \edge[10,0] \\  \\   \\
\place \edge[0,11] &&&&&&\dropup{3}{\langle 0\rangle}
&&&\dropup{3}{\langle Z_{12}\rangle} &&\place \\   \\
  &&&\dropleft{3}{\langle 0\rangle} &&&\plc &\plc &&\plc
&\plb \\
 &&&&&&\plc &\plc &&\plc &\plc \\   \\
 &&&\dropleft{3}{\langle Y_{21}\rangle} &&&\plc &\plc
&&\plc &\plb \\
 &&&&&&\plb &\plc &&\plb &\plc \\ 
 &&&&\place
}$$
\endignore
\medskip

We now conclude from Theorem 3.5 that we have found all the
$H$-primes of $A= \Oq(M_2(k))$. There are 14 in total, which we can
display as follows:

\ignore
$$\xymatrixrowsep{0.1pc}\xymatrixcolsep{0.1pc}
\xymatrix{
 &\plb &\plb &&&\hrzvrt &&\plc &\plb &&\plb &\plc &&&\plc
&\plc &&\plc &\plb \\
 &\plb &\plb &&&\hrz &&\plc &\plb &&\plb &\plc &&&\plc
&\plc &&\plc &\plc \\
\place &&&&& &&& &&& &&&& &&& &&\place \\
 && &&&\plc &\plc &&\plc &\plb &&\plb &\plc &&&\plc &\plc
&&\plc &\plb \\
 && &&&\plb &\plb &&\plb &\plb &&\plb &\plb &&&\plb &\plc
&&\plb &\plc \\   \\
 && &&&\plb &\plb &&\plb &\plb &&\plb &\plb  \\
 && &&&\plc &\plc &&\plc &\plb &&\plb &\plc
}$$
\endignore
\medskip

\noindent For a display showing the inclusions among these ideals, see
\cite{\GLmurcia, (3.6)}.
\enddefinition

\head 5. Appendix. Relations in $\OqMn$ \endhead

The proofs in this paper rely on a number of relations
among the generators and quantum minors in quantum matrix algebras. We
record and/or derive those relations in this appendix. Throughout, let
$A= \OqMn$ with
$k$ an arbitrary field and
$q\in \kx$ an arbitrary nonzero scalar.

\definition{5.1} {\bf (a)} We present the algebra $A$ with generators
$X_{ij}$ for $i,j=1,\dots,n$ and relations
$$\xalignat2 X_{ij}X_{lj} &= qX_{lj}X_{ij} &&\qquad (i<l)\\
X_{ij}X_{im} &= qX_{im}X_{ij} &&\qquad (j<m)\\
X_{ij}X_{lm} &= X_{lm}X_{ij} &&\qquad (i<l,\, j>m)\\
X_{ij}X_{lm} - X_{lm}X_{ij} &= (q - q^{-1})X_{im} X_{lj} &&\qquad
(i<l,\, j<m). \endxalignat$$ 
As is well known, $A$ is in fact a bialgebra, with
comultiplication $\Delta : A\rightarrow A\otimes A$ and counit
$\epsilon : A\rightarrow k$ such that 
$$\Delta(X_{ij})= \sum_{l=1}^n X_{il}\otimes X_{lj}
\qquad\text{and}\qquad \epsilon(X_{ij})= \delta_{ij}$$
for all $i,j$.

{\bf (b)} The algebra $A$ possesses various symmetries. In particular,
it supports a $k$-algebra automorphism $\bftau$ such that
$\bftau(X_{ij}) =X_{ji}$ for all $i,j$ \cite{\PaWa, Proposition
3.7.1}. Pairs of results which imply each other just through
applications of $\bftau$ will be simply referred to as ``symmetric''.

{\bf (c)} Given $U,V \subseteq \{1,\dots,n\}$ with $|U|=|V|=t$, let
$\OqMUV$ denote the $k$-subalgebra of $A$ generated by those $X_{ij}$
with $i\in U$ and $j\in V$. There is a natural isomorphism $\Oq(M_t(k))
\rightarrow \OqMUV$, which sends the quantum determinant of
$\Oq(M_t(k))$ to the quantum minor in $A$ involving rows from $U$ and
columns from $V$. As in \cite{\Duke}, we denote this quantum minor by
$[U\mid V]$, or in the form $[u_1\cdots u_t\mid v_1\cdots v_t]$ if we
wish to list the elements of $U$ and $V$. 

We shall also use the
isomorphism
$\Oq(M_t(k))
\rightarrow \OqMUV$ to simplify various proofs, since it will allow us
to work in smaller quantum matrix algebras than $A$ on occasion. For
example, since the quantum determinant in $\Oq(M_t(k))$ is central
(e.g., \cite{\PaWa, Theorem 4.6.1}), the quantum minor $[U\mid V]$
commutes with $X_{ij}$ for all $i\in U$ and $j\in V$. Consequently,
$$\xalignat2 [U\mid V][I\mid J] &= [I\mid J][U\mid V] &&\qquad
(I\subseteq U,\ J\subseteq V). \endxalignat$$

{\bf (d)} Recall from \cite{\NYM, Equation 1.9} the comultiplication
rule for quantum minors:
$$\Delta \bigl( [I\mid J] \bigr)= \sum_{|K|=|I|} [I\mid K] \otimes
[K\mid J].$$
\enddefinition

\definition{5.2} We next restate some
identities from
\cite{\PaWa}, given there for generators and maximal minors, in a form
that applies to minors of arbitrary size. Note the difference between
our choice of relations for $A$ (see (5.1)(a)) and that in
\cite{\PaWa, p\. 37}. Because of this, we must interchange $q$ and
$q^{-1}$ whenever carrying over a formula from \cite{\PaWa}.
\enddefinition

\proclaim{Lemma} Let $r,c\in
\{1,\dots,n\}$ and
$I,J\subseteq
\{1,\dots,n\}$ with $|I|=|J|\ge1$.

{\rm (a)} If $r\in I$ and $c\in J$, then $X_{rc}[I\mid J]= [I\mid
J]X_{rc}$.

{\rm (b)} If $r\in I$ and $c\notin J$, set $J^+= J\sqcup\{c\}$. Then
$$\allowdisplaybreaks
\align X_{rc}[I\mid J] -q^{-1}[I\mid J]X_{rc} &= (q^{-1}-q) \sum \Sb
j\in J\\ j>c\endSb (-q)^{-|J\cap [c,j]|} [I\mid 
J^+\setminus \{j\}] X_{rj} \tag1\\
[I\mid J]X_{rc} -qX_{rc}[I\mid J] &= (q-q^{-1}) \sum \Sb j\in
J\\ j>c\endSb (-q)^{|J\cap [c,j]|} X_{rj}[I\mid J^+\setminus \{j\}]
\tag2 \endalign$$

{\rm (c)} If $r\notin I$ and $c\in J$, set $I^+= I\sqcup\{r\}$. Then
$$\allowdisplaybreaks
\align X_{rc}[I\mid J] -q^{-1}[I\mid J]X_{rc} &= (q^{-1}-q) \sum \Sb
i\in I\\ i>r\endSb (-q)^{-|I\cap [r,i]|} [I^+\setminus \{i\}\mid
J]X_{ic} \tag1\\ 
[I\mid J]X_{rc} -qX_{rc}[I\mid J] &= (q-q^{-1}) \sum
\Sb i\in I\\ i>r\endSb (-q)^{|I\cap [r,i]|} X_{ic}[I^+\setminus
\{i\}\mid J] \tag2 \endalign$$
\endproclaim

\demo{Proof} Part (a) is clear. To obtain part (b), we may work in
$\Oq(M_{I\sqcup \{r'\},J^+}(k))$ for some $r'\notin
I$, and so there is no loss of generality in assuming that
$I\subset J^+=
\{1,\dots,n\}$. The desired relations then follow from the second
cases of parts (1) and (2) of \cite{\PaWa, Lemma 4.5.1}. Part (c) is
symmetric to (b).
\qed\enddemo

\proclaim{5.3\. Lemma} Let $U,V \subseteq
\{1,\dots,n\}$ with
$|U|=|V|$, and let $u_1,u_2\in U$ and $v_1,v_2\in V$. Set $U_s=
U\setminus \{u_s\}$ and $V_s= V\setminus \{v_s\}$ for $s=1,2$.

{\rm (a)} If $u_1<u_2$, then $[U_1\mid V_1] [U_2\mid V_1] =q^{-1}
[U_2\mid V_1] [U_1\mid V_1]$.

{\rm (b)} If $v_1<v_2$, then $[U_1\mid V_1] [U_1\mid V_2] =q^{-1}
[U_1\mid V_2] [U_1\mid V_1]$.

{\rm (c)} If $u_1<u_2$ and $v_1>v_2$, then $[U_1\mid V_1] [U_2\mid
V_2] = [U_2\mid V_2] [U_1\mid V_1]$.

{\rm (d)} If $u_1<u_2$ and $v_1<v_2$, then
$$[U_1\mid V_1] [U_2\mid V_2] -[U_2\mid V_2] [U_1\mid V_1] =
(q^{-1}-q) [U_2\mid V_1] [U_1\mid V_2].$$
\endproclaim

\demo{Proof} Since we may work in $\Oq(M_{U,V}(k))$, there is no loss
of generality in assuming that $U=V= \{1,\dots,n\}$. The result then
follows from \cite{\PaWa, Theorem 5.2.1}. \qed\enddemo

\definition{5.4} We also require the form of the {\it $q$-Laplace
relations\/} given in \cite{\NYM}. For index sets $I$ and $J$, set
$$\ell(I;J)= \bigl| \{(i,j)\in I\times J\mid i>j\} \bigr|.$$
\enddefinition

\proclaim{Lemma} {\rm ($q$-Laplace relations)} Let $I,J\subseteq
\{1,\dots,n\}$. 

{\rm (a)} If $J_1,J_2\subseteq \{1,\dots,n\}$ with $|J_1|+|J_2| =|I|$,
then
$$\sum \Sb I_1\sqcup I_2=I\\ |I_\nu|=|J_\nu|
\endSb (-q)^{\ell(I_1;I_2)} [I_1\mid J_1][I_2\mid J_2] =\cases
(-q)^{\ell(J_1;J_2)}[I\mid J_1\sqcup J_2] &\quad (J_1\cap J_2=
\varnothing)\\ 0 &\quad (J_1\cap J_2\ne \varnothing). \endcases$$

{\rm (b)} If $I_1,I_2\subseteq \{1,\dots,n\}$ with  $|I_1|+|I_2| =|J|$,
then
$$\sum \Sb J_1\sqcup J_2=J\\ |J_\nu|=|I_\nu|
\endSb (-q)^{\ell(J_1;J_2)} [I_1\mid J_1][I_2\mid J_2] =\cases
(-q)^{\ell(I_1;I_2)}[I_1\sqcup I_2\mid J] &\quad (I_1\cap I_2=
\varnothing)\\ 0 &\quad (I_1\cap I_2\ne \varnothing). \endcases$$
\endproclaim

\demo{Proof} \cite{\NYM, Proposition 1.1}.
\qed\enddemo

\definition{5.5} The $q$-Laplace relations simplify somewhat when one
of the index sets is a singleton, as follows.
\enddefinition

\proclaim{Corollary} Let $r,c\in \{1,\dots,n\}$ and
$I,J\subseteq
\{1,\dots,n\}$.

{\rm (a)} If $|I|= |J|+1$, then
$$\allowdisplaybreaks
\align \sum_{i\in I} (-q)^{|[1,i)\cap I|} X_{ic} [I\setminus \{i\}
\mid J] &= \cases (-q)^{|[1,c)\cap J|} [I\mid J\sqcup \{c\}] &\quad
(c\notin J)\\  0 &\quad (c\in J) \endcases \tag1\\
\sum_{i\in I} (-q)^{|(i,n]\cap I|} [I\setminus \{i\} \mid J] X_{ic} &=
\cases (-q)^{|(c,n]\cap J|} [I\mid J\sqcup \{c\}] &\quad (c\notin
J)\\  0 &\quad (c\in J). \endcases \tag2 \endalign$$

{\rm (b)} If $|J|= |I|+1$, then
$$\allowdisplaybreaks
\align \sum_{j\in J} (-q)^{|[1,j)\cap J|} X_{rj} [I\mid J\setminus
\{j\}] &= \cases (-q)^{|[1,r)\cap I|} [I\sqcup \{r\} \mid J] &\quad
(r\notin I)\\  0 &\quad (r\in I) \endcases \tag1\\
\sum_{j\in J} (-q)^{|(j,n]\cap J|} [I\mid J\setminus \{j\}] X_{rj} &=
\cases (-q)^{|(r,n]\cap I|} [I\sqcup \{r\} \mid J] &\quad (r\notin
I)\\  0 &\quad (r\in I). \endcases \tag2 \endalign$$
\endproclaim

\demo{Proof} (a) For the first case, fix $J_1= \{c\}$ and $J_2= J$. We
will use Lemma 5.4(a), which involves a sum over $I_1\sqcup I_2= I$
with $|I_1| =1$; thus $I_1= \{i\}$ and $I_2= I\setminus \{i\}$ for
some $i\in I$. In that case, $\ell(I_1;I_2)= |[1,i)\cap I|$ and
$\ell(J_1;J_2)= |[1,c)\cap J|$. Thus, formula (1) follows directly
from Lemma 5.4(a). Formula (2) follows similarly, where this time we
fix $J_1=J$ and $J_2= \{c\}$.

(b) These follow from (a) by symmetry. \qed \enddemo

\proclaim{5.6\. Lemma} Let $U,V \subseteq
\{1,\dots,n\}$ with
$|U|=|V|$.

{\rm (a)} Let $U= I\sqcup K$, and let $J_1,J_2
\subseteq V$ such that $|J_1| +|J_2|= 2|I| +|K|$. Then
$$\multline \sum \Sb K= K'\sqcup K''\\ |K'|= |J_1|-|I| \endSb
(-q)^{\ell(I;K')+
\ell(K'; K''\sqcup I)} [I\sqcup K' \mid J_1] [K''\sqcup I \mid J_2]\\
 =\cases (-q)^{ \ell(J_1\cap J_2; J_1\setminus J_2)+ \ell(J_1\setminus
J_2; J_2)} [I\mid  J_1\cap J_2] [U\mid V] &\quad (|J_1\cap J_2| =|I|)\\
 0 &\quad (|J_1\cap J_2| >|I|). \endcases
\endmultline$$

{\rm (b)} Let $V= J\sqcup L$, and let $I_1,I_2
\subseteq U$ such that $|I_1| + |I_2|= 2|J| +|L|$. Then
$$\multline \sum \Sb L= L'\sqcup L''\\ |L'|= |I_1|-|J| \endSb
(-q)^{\ell(J;L') + \ell(L';L''\sqcup J)} [I_1 \mid J \sqcup L'] [I_2
\mid J\sqcup L''] \\ 
 = \cases (-q)^{\ell(I_1 \cap I_2;I_1\setminus I_2) + \ell(I_1\setminus
I_2; I_2)} [I_1 \cap I_2\mid J] [U\mid V] &\quad (|I_1\cap I_2| =|J|)\\
 0 &\quad (|I_1\cap I_2| >|J|). \endcases
\endmultline$$ \endproclaim

\demo{Proof} By symmetry, we need only prove (a). Note that
$$|J_1| +|J_2|= 2|I| +|K|= |I|+|U| \ge |I|+|J_1\cup J_2|= |I|+|J_1|
+|J_2|-|J_1\cap J_2|,$$
whence $|J_1\cap J_2| \ge |I|$. Let (1) denote the left hand side of
the formula to be established, and (2) the first choice on the
right hand side.

 Expand each term $(-q)^{\ell(I;K')} [I\sqcup K' \mid J_1]$
using Lemma 5.4(b) and insert into (1). Thus, (1) equals
$$\sum \Sb K= K'\sqcup K''\\ J_1= J_1'\sqcup J_1'' \endSb (-q)^{
\ell(J_1'; J_1'')+ \ell(K'; K''\sqcup I)} [I\mid J_1'] [K'\mid J_1'']
[K''\sqcup I \mid J_2]. \tag3$$

We next claim that the sum
$$\sum \Sb U= I_1\sqcup I_2\\ J_1= J_1'\sqcup J_1'' \endSb
(-q)^{\ell(I_1;I_2)+ \ell(J_1'; J_1'')} [I\mid J_1'] [I_1\mid J_1'']
[I_2\mid J_2]  \tag4$$
equals (3). For $I_1$ as in (4), we have $|I|+|I_1|= |J'_1|+|J''_1|=
|J_1|$, and so Lemma 5.4(b) gives
$$\xalignat2 \sum_{J_1= J_1'\sqcup J_1''} (-q)^{\ell(J_1'; J_1'')}
[I\mid J_1'] [I_1\mid  J_1''] &=0 &&\qquad (I_1 \not\subseteq K),
 \tag5 \endxalignat$$
because $I\cap I_1 \ne \varnothing$ in this case. On the other hand,
for fixed $J'_1$, $J''_1$ as in (4), we have
$$\multline \sum \Sb U= I_1\sqcup I_2\\ I_1\subseteq K \endSb
(-q)^{\ell(I_1;I_2)} [I_1\mid J_1''] [I_2\mid J_2]\\
 = \sum_{K= K'\sqcup K''} (-q)^{\ell(K'; K''\sqcup I)} [K'\mid J_1'']
[K''\sqcup I \mid J_2]. \endmultline \tag6$$
It follows from (5) and (6) that $(4)=(3)$ as claimed, and thus
$(1)=(4)$.

For $J_1''$ as in (4), we have $|J_1''|= |K'|= |K|- |K''|= |U|- |K''
\sqcup I|= |V|- |J_2|$. Hence, Lemma 5.4(a) says that 
$$\sum_{ U= I_1\sqcup I_2} (-q)^{\ell(I_1;I_2)} [I_1\mid
J_1''] [I_2\mid J_2]= \cases (-q)^{\ell(J_1'';J_2)} [U\mid V] &\quad
(J_1''\cap J_2 =\varnothing)\\ 0 &\quad (J_1''\cap J_2 \ne \varnothing).
\endcases  \tag7$$
Substituting (7) into (4), it follows that (1) is equal to the sum
$$\sum_{J_1= J_1'\sqcup J_1''} (-q)^{\ell(J_1'; J_1'')+ \ell(J_1'';J_2)}
[I\mid J_1'] d(J_1''),  \tag8$$
where $d(J_1'')= [U\mid V]$ if $J_1''$ and $J_2$ are disjoint, but
$d(J_1'')= 0$ otherwise.

If $|J_1\cap J_2| >|I|$, then since any $|J_1'|= |I|$, we see that
$J_1\cap J_2 \nsubseteq J_1'$ and so $J_1''\cap J_2 \ne\varnothing$.
Thus in this case all $d(J_1'')=0$, and so $(1)=(8)=0$.

Finally, suppose that $|J_1\cap J_2| =|I|$. Then the only time $J_1''$
and $J_2$ can be disjoint is when $J_1'= J_1\cap J_2$, and therefore
$(1)= (8)= (2)$ in this case. \qed\enddemo

\proclaim{5.7\. Lemma} Let $r,c\in
\{1,\dots,n\}$ and
$I,J\subseteq
\{1,\dots,n\}$ with $|I|=|J|\ge1$. If $r> \max(I)$ and $c> \max(J)$,
then
$$[I\mid J]X_{rc}- q^2 X_{rc}[I\mid J]= (1-q^2)[I\sqcup \{r\} \mid
J\sqcup \{c\}].$$
\endproclaim

\demo{Proof} Since we may work in $\Oq(M_{I\sqcup\{r\},
J\sqcup\{c\}}(k))$, it suffices to consider the case that
 $r=c=n$ and $I=J=
\{1,\dots, n-1\}$. Now $[I\mid J]= A(n\,n)$ in the
notation of \cite{\PaWa, (4.3)}. Set $D_q= [I\sqcup \{r\} \mid
J\sqcup \{c\}]= [1\cdots n\mid 1\cdots n]$.

The first two $q$-Laplace relations in \cite{\PaWa, Corollary 4.4.4}
yield
$$\sum_{j=1}^n (-q)^{j-n} X_{nj} A(n\,j)= \sum_{j=1}^n (-q)^{n-j}
A(n\,j) X_{nj} =D_q. \tag1$$
Solving for $X_{nn}A(n\,n)$ and $A(n\,n)X_{nn}$, we obtain
$$\align X_{nn}A(n\,n) &= D_q- \sum_{j<n} (-q)^{j-n} X_{nj} A(n\,j)
\tag2 \\
A(n\,n)X_{nn} &= D_q- \sum_{j<n} (-q)^{n-j}
A(n\,j) X_{nj}. \tag3 \endalign$$

For any $j$, the third relation of \cite{\PaWa, Lemma 5.1.2} implies
that
$$X_{nj}A(n\,j)= A(n\,j) X_{nj}+ (1-q^{-2}) \sum_{l<j} (-q)^{j-l}
A(n\,l) X_{nl}. \tag4$$
Substituting (4) into (2) for all $j<n$, we obtain
$$\aligned X_{nn} A(n\,n) &= D_q- \sum_{j<n} (-q)^{j-n} A(n\,j)
X_{nj}\\
 &\qquad -(1-q^{-2}) \sum_{j<n} \sum_{l<j} (-q)^{2j-l-n} A(n\,l)
X_{nl}\\
 &= D_q- \sum_{l<n} \biggl[ (-q)^{l-n}+ (1-q^{-2}) \sum_{l<j<n}
(-q)^{2j-l-n} \biggr] A(n\,l) X_{nl}. \endaligned \tag5$$
The expression in square brackets can be simplified as follows:
$$\aligned (-q)^{l-n}+ (1-q^{-2}) \sum_{l<j<n} (-q)^{2j-l-n} &=
(-q)^{l-n} \biggl[ 1+ (1-q^{-2})\sum_{0<m<n-l} (-q)^{2m} \biggr]\\
 &= (-q)^{n-l-2}. \endaligned \tag6$$ 
Substituting (6) into (5) and replacing $l$ by $j$, we obtain
$$X_{nn} A(n\,n)= D_q- \sum_{j<n} (-q)^{n-j-2} A(n\,j) X_{nj}. \tag7$$
Finally, combining (3) with (7), we conclude that
$$A(n\,n) X_{nn}- q^2 X_{nn} A(n\,n)= (1-q^2) D_q,$$
as desired. \qed\enddemo

\enddocument